\pdfoutput=1
\RequirePackage{ifpdf}
\ifpdf
\documentclass[pdftex]{sigma}
\else
\documentclass{sigma}
\fi

\usepackage[all]{xy}

\numberwithin{equation}{section}

\newtheorem{Theorem}{Theorem}[section]
\newtheorem{Corollary}[Theorem]{Corollary}
\newtheorem{Lemma}[Theorem]{Lemma}
\newtheorem{Proposition}[Theorem]{Proposition}
{\theoremstyle{definition}
\newtheorem{Definition}[Theorem]{Definition}
\newtheorem{Remark}[Theorem]{Remark}
\newtheorem{Example}[Theorem]{Example}
}

\begin{document}

\allowdisplaybreaks

\renewcommand{\thefootnote}{$\star$}

\renewcommand{\PaperNumber}{018}

\FirstPageHeading

\ShortArticleName{Fukaya Categories as Categorical Morse Homology}

\ArticleName{Fukaya Categories as Categorical Morse Homology\footnote{This paper is a~contribution to the Special Issue ``Mirror
Symmetry and Related Topics''.
The full collection is available at
\href{http://www.emis.de/journals/SIGMA/mirror_symmetry.html}{http://www.emis.de/journals/SIGMA/mirror\_{}symmetry.html}}}

\Author{David NADLER}

\AuthorNameForHeading{D.~Nadler}

\Address{Department of Mathematics, University of California, Berkeley, Berkeley, CA 94720-3840, USA}
\Email{\href{mailto:nadler@math.berkeley.edu}{nadler@math.berkeley.edu}}
\URLaddress{\url{http://math.berkeley.edu/~nadler/}}

\ArticleDates{Received May 16, 2012, in f\/inal form February 21, 2014; Published online March 01, 2014}

\Abstract{The Fukaya category of a~Weinstein manifold is an intricate symplectic inva\-riant of high interest in mirror symmetry
and geometric representation theory.
This paper informally sketches how, in analogy with Morse homology, the Fukaya category might result from gluing together Fukaya
categories of Weinstein cells.
This can be formalized by a~re\-collement pattern for Lagrangian branes parallel to that for constructible sheaves.
Assuming this structure, we exhibit the Fukaya category as the global sections of a~sheaf on the conic topology of the Weinstein
manifold.
This can be viewed as a~symplectic analogue of the well-known algebraic and topological theories of (micro)localization.}


\Classification{53D37}

\renewcommand{\thefootnote}{\arabic{footnote}}
\setcounter{footnote}{0}

\vspace{-2mm}

\section{Introduction}

To realize ``compact, smooth'' global objects as glued together from simpler local pieces, one often pays the price that the
local pieces are ``noncompact'' or ``singular''.
For several representative examples, one could think about compact manifolds versus cells and simplices, smooth projective
varieties versus smooth af\/f\/ine varieties and singular hyperplane sections, vector bundles with f\/lat connection versus regular
holonomic ${\mathcal D}$-modules, or perhaps most universally of all, irreducible modules versus induced modules.

In this paper, we informally sketch how a~similar pattern might hold for exact Lagrangian branes in Weinstein manifolds: it
should be possible to systematically glue together (in the homological sense) complicated branes from simple but noncompact
branes.
If the initial branes are compact, then Seidel's theory of iterated Lefschetz f\/ibrations and their vanishing cycles and
thimbles~\cite{Seidel} provides a~complete solution, both theoretically and computationally.
For applications to some questions in geometric representation theory, it is useful to go further and consider noncompact branes
from the start.
When the Weinstein manifold is a~cotangent bundle, constructible sheaves capture the structure of all branes, compact and
noncompact alike, and in particular, organize the gluing relations between local and global calculations~\cite{ N, NSpr, NZ}.
We turn here to a~basic structure of a~Weinstein manifold and discuss decomposing branes along its unstable coisotropic Morse
cells.
One can interpret the resulting proposed presentation of the Fukaya category as a~categorical form of Morse homology.

As a~consequence, assuming such a~gluing pattern, we derive a~canonical localization of branes: we construct a~sheaf of
categories whose global sections recover the Fukaya category.
This provides a~symplectic counterpart to the algebraic theory of Beilinson--Bernstein localization~\cite{BB} and the topological
microlocalization of Kashiwara--Schapira~\cite{KS}.
Though we are exclusively occupied here with general structure, we anticipate applications of potentially broad appeal.
The localization of branes should assist in the calculation of many Fukaya categories, in particular in terms of algebraic
quantizations when the Weinstein manifold is the Hamiltonian reduction of a~cotangent bundle.
For the simplest example, one could consider Slodowy slices and their smoothings: on the one hand, their compact branes have
been studied extensively in the context of knot homologies~\cite{KhS, M, SS}; on the other hand, their noncommutative modules
globalize to modules over Lie algebras and f\/inite $W$-algebras (see~\cite{kost, Lo1, Lo2, premet} for origins and up-to-date
lists of references).
More broadly, one could also consider other symplectic settings for geometric representation theory, for example for
localizations of Cherednik-type algebras and hypertoric enveloping algebras~\cite{BLPW,GS1, GS2, KRouquier, MN}.
This paper suggests useful tools to see that such targets admit a~single categorical quantization which can be described in
terms of localized branes or alternatively noncommutative modules.

From a~concrete geometric perspective, we have been particularly inspired by Seidel's theory of exact symplectic Dehn
twists~\cite{S01b,S01, S03, Seidel}.
For our understanding of the foundations of the subject, we have turned to the work of Eliashberg, Gromov, and
Weinstein~\cite{E, EG, W}, Fukaya, Oh, Ohta, and Ono~\cite{FOOO}, Seidel~\cite{Seidel}, Fukaya and Oh~\cite{FO}, Wehrheim and
Woodward~\cite{ww6,ww3, ww1,ww4, ww2, ww5}, Kashiwara and Schapira~\cite{KS} and Lurie~\cite{HA, topos, TFT}.
There is a~wealth of results in closely aligned directions, in particular in the analysis of the homology of closed strings.
We are familiar with only a~small part of this growing literature, in particular, the guiding advances of Seidel~\cite{S08,
S09}, and the striking results of Bourgeois, Cieliebak, Eckholm, and Eliashberg~\cite{BEE2, BEE1, C}.

\looseness=-1
This paper outlines anticipated structure and does not provide comprehensive arguments.
Due to the high technical demands, our discussion of functors on noncompact branes is only a~sketch highlighting our main
expectations within this rapidly developing subject.
Assu\-ming these structural results, we provide complete arguments for the localization of the Fukaya cate\-gory.

In what follows, we f\/irst outline the specif\/ic setup adopted in this paper, then go on to describe our main expectations,
results and their immediate precursors.

\subsection{Setup} Let $(M, \theta)$ be a~Weinstein manifold, so a~manifold $M$ equipped with a~one-form $\theta$ whose
dif\/ferential $\omega = d\theta$ is a~symplectic form such that there exists a~suitably compatible Morse function $h:M\to \mathbb
R$.
(We refer the reader to Section~\ref{sect weinstein} below for more details on the compatibility, but note here the potential
point of confusion: we prefer the f\/lexibility of not including the choice of a~particular Morse function in the data, though it
is standard to do so in some other contexts.) We will always assume that $M$ is real analytic, and all subsets and functions are
subanalytic (or def\/inable within some f\/ixed o-minimal context).
The basic source of examples are Stein manifolds, or more specif\/ically, smooth af\/f\/ine complex varieties.

Let $Z$ be the Liouville vector f\/ield on $M$ characterized by $\theta = i_Z\omega$, and let $\mathfrak c\subset M$ be the f\/inite
subset of zeros of $Z$.
For generic data, the f\/low of $Z$ provides a~stratif\/ication ${\mathcal S} = \{C_p \}_{p\in \mathfrak c}$ into coisotropic
unstable cells $C_p\subset M$ contracting to the zeros $p\in \mathfrak c$.
Hamiltonian reduction along each coisotropic unstable cell $C_p\subset M$ produces a~contractible Weinstein manifold $(M_p,
\theta_p)$ which we refer to as a~Weinstein cell.

Let $F(M)$ denote the Fukaya category of not necessarily compact exact Lagrangian branes, and let $\operatorname{Perf} F(M)$
denote the stable Fukaya category of perfect modules over $F(M)$.
Here we f\/ix a~coef\/f\/icient f\/ield $k$, and by perfect modules, mean summands of f\/inite complexes of representable modules with
values in $k$-chain complexes.
Let us brief\/ly orient the reader as to what version of the Fukaya category we work with.
(See Section~\ref{sec fukaya} below for a~more detailed discussion of all of the following notions.)

First, recall that any symplectic manifold $M$ is canonically oriented, and one can unambiguously speak about its Chern classes.
To work with graded Lagrangian submanifolds, we will make the standard assumption that the characteristic class $2 c_1(M)$ is
trivialized.
For simplicity, we will also assume that $M$ comes equipped with a~spin structure.
Then by Lagrangian brane, we will mean a~graded Lagrangian submanifold equipped with a~f\/inite-dimensional local system and pin
structure.

Second, by the Fukaya category $F(M)$, we mean the inf\/initesimal (as describes the setting of~\cite{N, NZ, S01, Seidel}) rather
than wrapped variant (as found in~\cite{AS}; this is one extreme of the partially wrapped paradigm~\cite{AurICM, Aur}).
Our perturbation framework involves Hamiltonian isotopies (in the direction of the rotated Liouville vector f\/ield) of constant
size with respect to a~radial coordinate near inf\/inity rather than of linear growth.
If one specif\/ies a~conic Lagrangian support $\Lambda \subset M$ for branes, there is a~resulting full subcategory
$\operatorname{Perf}_\Lambda F(M) \subset \operatorname{Perf} F(M)$.
The corresponding partially wrapped category might be viewed as its ``Verdier dual'' (the relation appears analogous to that of
cohomology and homology, or more immediately, perfect and coherent ${\mathcal O}$-modules).
In fact, one might expect the partially wrapped category to embed into the category of modules over $\operatorname{Perf}_\Lambda
F(M)$.
Thus by coupling results for the inf\/initesimal category with notions of support and homological bounds, one might expect to
obtain parallel results for the partially wrapped category.
We have included some further discussion of expectations in Remark~\ref{wrapped} below.

\subsection{Recollement for branes}

Our f\/irst goal is to outline how the stable Fukaya category $\operatorname{Perf} F(M)$ might be recoverable from the stable
Fukaya categories $\operatorname{Perf} F(M_p)$ of its Weinstein cells together with gluing data in the form of natural
adjunctions.

Given a~closed coisotropic submanifold $i:C\to M$ that is a~union of unstable cells, we will sketch an expected semiorthogonal
decomposition of $\operatorname{Perf} F(M)$.
On the one hand, the formalism of Hamiltonian reduction provides a~Lagrangian correspondence
\begin{gather*}
\xymatrix{ N & \ar@{->>}[l]_-{q} C \ar@{^(->}[r]^-{i} & M, }
\end{gather*}
where $N$ is a~Weinstein manifold, and $q$ is the quotient along the integrable isotropic foliation determined by $i$.
On the other hand (as explained in Lemma~\ref{lemma open restriction} below), the open complement $j:M^\circ = M \setminus C\to
M$ naturally inherits the structure of a~Weinstein manifold.
(Note that we do not include a~specif\/ic compatible Morse function in the data of a~Weinstein manifold.
The compatible Morse function we construct for $M^\circ$ is not simply the restriction of a~compatible Morse function for $M$.)

We will outline how there might be a~recollement pattern for Lagrangian branes analogous to standard gluings for constructible
sheaves (as recalled in Section~\ref{dev} below): it is a~natural diagram of adjunctions
\begin{gather}
\label{eq functors}
\xymatrix{ \ar[rr]^-{\mathfrak i_! \simeq \mathfrak i_*} \operatorname{Perf} F(N) && \operatorname{Perf} F(M)
\ar@/^2pc/[ll]_-{i^!}\ar@/_2pc/[ll]_-{\mathfrak i^*} \ar[rr]^-{\mathfrak j^! \simeq \mathfrak j^*} && \operatorname{Perf}
F(M^\circ) \ar@/^2pc/[ll]_-{\mathfrak j_*} \ar@/_2pc/[ll]_-{\mathfrak j_!} }
\end{gather}
with $\mathfrak i_! \simeq \mathfrak i_*$, $\mathfrak j_!$, $\mathfrak j_*$ fully faithful embeddings, along with exact triangles
of functors
\begin{gather}
\label{eq triangles}
\xymatrix{ \mathfrak i_{!} \mathfrak i^! \ar[r]^-c & \operatorname{id}_{\operatorname{Perf} F(M)} \ar[r]^-u & \mathfrak
j_*\mathfrak j^* \ar[r]^-{[1]} &, & \mathfrak j_{!}\mathfrak j^! \ar[r]^-c & \operatorname{id}_{\operatorname{Perf} F(M)}
\ar[r]^-u & \mathfrak i_*\mathfrak i^* \ar[r]^-{[1]}& },
\end{gather}
where $c$ denotes the counits of adjunctions and $u$ the units.

From a~symplectic perspective, the functors involved are ``Dehn twists around cells''.
Alternatively, one can also interpret them as resulting from Lagrangian surgery.

\begin{Remark}
If $M$ admits an anti-symplectomorphism compatible with other structures, then one can obtain an analogue of Verdier duality as
well.
\end{Remark}

\begin{Remark}
\label{mv}
The exact triangles are categorical versions of the long exact sequences of cohomology of a~pair.
It is useful to observe that they are formally equivalent to exact triangles
\begin{gather*}
\xymatrix{ \operatorname{id}_{F(M)}\ar[r] & \mathfrak j_* \mathfrak j^* \oplus \mathfrak i_*\mathfrak i^* \ar[r] & \mathfrak
i_*\mathfrak i^* \mathfrak j_* \mathfrak j^* \ar[r]^-{[1]} &, & \mathfrak i_!i^!\mathfrak j_!\mathfrak j^! \ar[r] & \mathfrak
j_!\mathfrak j^! \oplus\mathfrak i_!i^! \ar[r] & \operatorname{id}_{F(M)}\ar[r]^-{[1]} & },
\end{gather*}
where the f\/irst and second maps are the units and counits of adjunctions.
These are categorical versions of the Mayer--Vietoris sequences of cohomology of a~covering pair.
\end{Remark}

\begin{Remark}
We mention here an attractive way to potentially repackage the recollement pattern.

Consider the dif\/ferential graded derived category $\operatorname{Sh}_{\mathcal S}(M)$ of complexes of sheaves constructible
along the stratif\/ication ${\mathcal S} = \{C_p \}_{p\in \mathfrak c}$ by coisotropic unstable cells $C_p \subset M$.
We equip $\operatorname{Sh}_{\mathcal S}(M)$ with its natural symmetric monoidal structure given by tensor product.
The monoidal unit is the constant sheaf $k_M \in \operatorname{Sh}_{\mathcal S}(M)$.
More generally, given a~union of unstable coisotropic cells $i:C\to M$, the extension by zero $k_{C!} = i_{!}k_{C} \in
\operatorname{Sh}_{{\mathcal S}}(M)$ is an idempotent.

Now given the recollement functors, it seems likely one should be able to reorganize them into a~natural fully faithful monoidal
embedding
\begin{gather*}
\xymatrix{ \operatorname{Sh}_{{\mathcal S}}(M) \ar[r] & \operatorname{End}_{{\rm st}_k} (\operatorname{Perf} F(M)) }
\end{gather*}
characterized by the property that $ k_{C!} \longmapsto \mathfrak i_!\mathfrak i^*.
$ This would allow one to analyze $\operatorname{Perf} F(M)$ as a~module over the ``spectrum'' of the commutative algebra
$\operatorname{Sh}_{\mathcal S}(M)$.
\end{Remark}

Such a~recollement pattern would provide various frameworks for gluing together $\operatorname{Perf} F(M)$ from the constituent
pieces $\operatorname{Perf} F(N)$, $\operatorname{Perf} F(M^\circ)$.
To pursue gluing, it is useful to regard each of the above categories as a~small stable idempotent-complete $k$-linear
$\infty$-category.
Then we can adopt the foundations of~\cite{HA, topos}, and work within the $\infty$-category ${\rm st}_k$ of small stable
idempotent-complete $k$-linear $\infty$-categories.

According to the reinterpretation of Remark~\ref{mv}, we can view $\operatorname{Perf} F(M)$ as classifying triples of data
$L^\circ \in \operatorname{Perf} F(M^\circ)$, $L_N \in \operatorname{Perf} F(N)$ together with a~morphism
\begin{gather*}
r \in \operatorname {Hom}_{\operatorname{Perf} F(N)}(\mathfrak i^! \mathfrak j_!L^\circ, L_N).
\end{gather*}
Let us recast this in monadic terms by considering the adjunction
\begin{gather*}
\xymatrix{ L =\mathfrak j_! \oplus\mathfrak i_!: \operatorname{Perf} F(M^\circ) \oplus \operatorname{Perf} F(N)\ar@<0.5ex>[r]
&\ar@<0.5ex>[l] \operatorname{Perf} F(M):R =\mathfrak j^! \oplus \mathfrak i^! }.
\end{gather*}
We obtain a~resulting monad, or in other words, algebra object in endomorphisms
\begin{gather*}
T = RL\in \operatorname{End}(\operatorname{Perf} F(M^\circ ) \oplus \operatorname{Perf} F(N)).
\end{gather*}
The $\infty$-categorical Barr--Beck theorem provides a~canonical equivalence
\begin{gather*}
\operatorname{Perf} F(M) \simeq \operatorname{Mod} _T(\operatorname{Perf} F(M^\circ ) \oplus \operatorname{Perf} F(N)),
\end{gather*}
where the right hand side denotes module objects over the monad $T$.

We can inductively apply the above considerations by successively taking the coisotropic submanifold $i:C\to M$ to be a~single
closed unstable cell.
Recall that $\mathfrak c\subset M$ denotes the f\/inite subset of zeros of the Liouville vector f\/ield $Z$.
For each $p\in \mathfrak c$, we have the Hamiltonian reduction of the corresponding unstable cell
\begin{gather*}
\xymatrix{ M_p & \ar@{->>}[l]_-{q_p} C_p \ar@{^(->}[r]^-{i_p} & M }.
\end{gather*}
By induction, the recollement pattern provides an adjunction
\begin{gather*}
\xymatrix{ L = \oplus_{p\in \mathfrak c} \mathfrak i_{p!}: \oplus_{p\in \mathfrak c} \operatorname{Perf} F(M_p) \ar@<0.5ex>[r]
&\ar@<0.5ex>[l] \operatorname{Perf} F(M):R = \oplus_{p\in \mathfrak c} \mathfrak i_p^! }.
\end{gather*}
We obtain a~resulting upper-triangular monad
\begin{gather*}
T = RL\in \operatorname{End}\big(\oplus_{p\in \mathfrak c} \operatorname{Perf} F(M_p)\big)
\end{gather*}
with matrix entries $\mathfrak i_q^! \mathfrak i_{p!}$ for unstable cells $C_q \subset \overline{C}_p$.
Finally, the $\infty$-categorical Barr--Beck theorem provides a~canonical equivalence
\begin{gather*}
\xymatrix{ \operatorname{Perf} F(M) \simeq \operatorname{Mod} _T(\oplus_{p\in \mathfrak c} \operatorname{Perf} F(M_p)) }.
\end{gather*}

One can view this presentation of $\operatorname{Perf} F(M)$ as a~categorif\/ied version of the Morse homology of $M$.
First, to each critical point $p\in \mathfrak c$, we assign the stable Fukaya category $\operatorname{Perf} F(M_p)$ of the
Weinstein cell $M_p$.
From the perspective of Morse homology, one can view $F(M_p)$ as the analogue of a~scalar vector space $k[\deg p]$ shifted by
the index of $p$.
(Since the $\infty$-category ${\rm st}_k$ of small stable idempotent-complete $k$-linear $\infty$-categories is not itself stable, it
is unsurprising that ``suspension'' in the current setting is not invertible as it is in the traditional setting.)

Second, to pairs of critical points $p, q\in \mathfrak c$, we assign the individual term $\mathfrak i_q^! \mathfrak i_{p!}$ of
the gluing monad $T$.
More generally, to collections of critical points, we assign the corresponding terms of the gluing monad $T$ together with the
monadic structure maps among them.
From the perspective of Morse homology, one can view the individual terms and their monadic structure as the analogue of
boundary maps and their higher relations.

Part of the appeal of Morse homology is that the boundary maps and their higher relations localize along the spaces of f\/low
lines connecting the relevant critical points.
Furthermore, the boundary maps and their higher relations admit simple descriptions in terms of the spaces of f\/low lines.
For instance, the boundary maps themselves are given by counts of isolated f\/low lines.
We will develop an abstract analogue of this picture directly below.

\subsection{Localization of branes} Now let us assume the above recollement pattern and see what it further implies.
We will deduce a~localization of the stable Fukaya category $\operatorname{Perf} F(M)$ over the conic topology of $M$.
There are many inspiring precedents~\cite{Ab, AS, SMorse, Spants, STZ} and informed assertions~\cite{kont} that suggest such
a~construction should be possible.
In our own thinking, we have often returned to the setting of cotangent bundles as a~guide: combining the equivalence of
Lagrangian branes with constructible sheaves~\cite{ N, NSpr, NZ} and the formalism of microlocalization~\cite{KS} leads to
a~complete solution in that case.

To describe the localization for a~general Weinstein manifold $(M, \theta)$, let us return to the geometry of the Liouville
vector f\/ield $Z$.
We say a~subset of $M$ is {\em conic} if it is invariant under the f\/low of $Z$, and use the term {\em core} to refer to the
compact conic isotropic subvariety $K\subset M$ of points that do not escape to inf\/inity under the f\/low.
(Many authors prefer the term ``skeleton'' for $K\subset M$, but we have opted for core so as not to confuse it with the
characteristic cone $\Lambda \subset M$ introduced below.
There is a~growing number of authors who would refer to $\Lambda \subset M$ as a~noncompact skeleton, and so we prefer to
distinguish the compact skeleton $K\subset M$ with the name core.) We use the term {\em ether} to refer to the complement $E =
M\setminus K $, and def\/ine the {\em projectivization} $M^\infty = E/{\mathbb R}_+$ to be the compact contact manifold obtained
by quotienting the ether by the f\/low.
More generally, given any conic subset $A\subset M$, we can consider its projectivization $A^\infty = (A\cap E)/{\mathbb
R}_+\subset M^\infty.$

To any object $L\in \operatorname{Perf} F(M)$, we assign its singular support $\operatorname{ss}(L) \subset M$ which is a~closed conic
isotropic subvariety depending only on the isomorphism class of the object (see Section~\ref{sing supp} below).
The singular support records the coarse homological nontriviality of objects, and for a~Lagrangian brane is a~subset of its
limiting dilation to a~conic Lagrangian subvariety.
For example, if $L$ is a~compact Lagrangian brane, then $\operatorname{ss}(L) \subset K$.

Let us f\/ix a~(most likely singular and noncompact) conic Lagrangian subvariety $\Lambda\subset M$ which we refer to as the {\em
characteristic cone}.
We will assume that $\Lambda$ contains the core $K \subset M$, and hence is completely determined by its projectivization
$\Lambda^\infty \subset M^\infty$.
Note that the inclusions $K\subset \Lambda \subset M$ are all homotopy equivalences, so what will interest us most is their
local geometry.

Let us consider the full subcategory $\operatorname{Perf}_\Lambda (M) \subset \operatorname{Perf} F(M)$ of objects $L\in
\operatorname{Perf} F(M)$ with singular support satisfying $\operatorname{ss}(L) \subset \Lambda$.
So for example, if $L$ is a~compact Lagrangian brane, then it provides an object of $\operatorname{Perf}_\Lambda F(M)$
irrespective of the choice of $\Lambda$.
Note that $\Lambda_1 \subset \Lambda_2$ implies $\operatorname{Perf}_{\Lambda_1} F(M) \subset \operatorname{Perf}_{\Lambda_2}
F(M)$, and $\operatorname{Perf} F(M) = \cup_{\Lambda} \operatorname{Perf}_\Lambda F(M)$.

Recall that $\mathfrak c\subset M$ denotes the f\/inite subset of zeros of the Liouville vector f\/ield $Z$.
For each $p\in \mathfrak c$, we have the Hamiltonian reduction of the corresponding unstable cell
\begin{gather*}
\xymatrix{ M_p & \ar@{->>}[l]_-{q_p} C_p \ar@{^(->}[r]^-{i_p} & M }.
\end{gather*}
We obtain a~conic Lagrangian subvariety $\Lambda_p\subset M_p$ by setting $\Lambda_p = q_p(i_p^{-1}(\Lambda))$, and similarly
a~full subcategory $\operatorname{Perf}_{\Lambda_p} F(M_p) \subset \operatorname{Perf} F(M_p)$.

Now by localizing $\operatorname{Perf}_\Lambda F(M)$ with respect to singular support, we obtain the following.
Its verif\/ication appeals to the recollement pattern outlined above.

\begin{Theorem}
\label{introtheorem2}
Assume the recollement pattern of diagrams~\eqref{eq functors} and~\eqref{eq triangles}.

There exists a~$\operatorname{st}_k$-valued sheaf ${\mathcal F}_\Lambda$ on the conic topology of $M$ with the following properties:
\begin{enumerate}\itemsep=0pt
\item[$(1)$] The support of ${\mathcal F}_\Lambda$ is the characteristic cone $\Lambda \subset M$.

\item[$(2)$] The global sections of ${\mathcal F}_\Lambda$ are canonically equivalent to $\operatorname{Perf}_\Lambda F(M)$.

\item[$(3)$] The restriction of ${\mathcal F}_\Lambda$ to an open Weinstein submanifold $M^\circ \subset M$ is canonically equivalent to
the sheaf ${\mathcal F}_{\Lambda^\circ}$ constructed with respect to $\Lambda^\circ = \Lambda \cap M^\circ$.

\item[$(4)$] For each zero $p\in \mathfrak c$, the sections of ${\mathcal F}_\Lambda$ lying strictly above the unstable cell $C_p \subset
M$ are canonically equivalent to $\operatorname{Perf}_{\Lambda_p} F(M_p)$.
\end{enumerate}
\end{Theorem}

\begin{Example}
Here is a~description of the sheaf ${\mathcal F}_\Lambda$ in the simplest example.

Consider the two-dimensional Weinstein cell $M = \mathbb C$ with standard Liouville form $\theta$ and projectivization $M^\infty
\simeq S^1$.
Its core is the single point $K = \{0\} \subset \mathbb C$, and its ether is the complement $E = \mathbb C^* \subset \mathbb C$.
Any characteristic cone $\Lambda \subset \mathbb C$ will be the union of $K = \{0\}$ with f\/initely many rays.
For $n = 0, 1, 2, \ldots$, let $\Lambda_n \subset \mathbb C$ denote the characteristic cone with $n$ rays.

Then for $n>1$, $\operatorname{Perf}_{\Lambda_n} F(M)$ is equivalent to f\/inite-dimensional modules over the $A_{n-1}$-quiver
$1\to 2 \to \dots \to n$.
For $n=0$ or $1$, it is the zero category.
This is also the stalk of the sheaf ${\mathcal F}_{\Lambda_n}$ at the point $0\in \mathbb C$.
Its stalk at other points $x\in \mathbb C$ is (not necessarily canonically) equivalent to $\operatorname{Perf} k$ when $x\in
\Lambda_n$, and is the zero category otherwise.

\end{Example}

\begin{Example}
As a~continuation of the previous example (or in fact a~direct generalization), we could take $M$ to be an open Riemann surface
of genus $g$ with $k>0$ punctures.
Then its projectivization $M^\infty$ is the disjoint union of $k$ circles, and its core $K\subset M$ is a~graph homotopy
equivalent to a~bouquet of $2g + (k-1)$ circles.
Any characteristic cone $\Lambda\subset M$ will be the union of $K$ with f\/initely many rays shooting from the nodes of $K$ of\/f
to the boundary.
So all together, we can view $\Lambda$ as a~``ribbon graph'': an abstract graph with some noncompact edges embedded in the
Riemann surface $M$ as a~retract.

Now the sheaf ${\mathcal F}_\Lambda$ is supported along $\Lambda$ (so equivalently can be viewed as a~sheaf on $\Lambda$).
It is locally constant along the edges of $\Lambda$ with stalks (not necessarily canonically) equivalent to $\operatorname{Perf}k$.
At the nodes of $\Lambda$, it is (not necessarily canonically) equivalent to f\/inite-dimensional modules over the $A_{n}$-quiver
where $n+2$ is the valency of a~node.
(So if a~node has valency $0$ or~$1$, the stalk is the zero category.) The global sections $\operatorname{Perf}_\Lambda(M)$ can
be easily calculated as the global sections of ${\mathcal F}_\Lambda$, with the combinatorial form of the answer dependent on
the presentation of the data of $\Lambda$ and its embedding in $M$.
(For example, this is compatible with the ``Òconstructible plumbing model'' for ribbon graphs with valency ${\leq}4$ developed
in~\cite{STZ}.)
\end{Example}

\begin{Remark}
It is possible to say much more about ${\mathcal F}_\Lambda$.
Let us content ourselves here with mentioning that when $M$ is a~cotangent bundle, sections of ${\mathcal F}_\Lambda$ are
equivalent to (the sheaf\/if\/ication of) microlocal sheaves.
In general, a~similar description holds for the restriction of ${\mathcal F}_\Lambda$ to the complement $M\setminus \mathfrak c$
of the zeros of the Liouville vector f\/ield $Z$.
Namely, the quotient $(M\setminus \mathfrak c)/{\mathbb R}_+$ by the Liouville f\/low is naturally a~(non-Hausdorf\/f) contact
manifold, and the restriction of ${\mathcal F}_\Lambda$ to the complement $M\setminus \mathfrak c$ can be obtained by pulling
back microlocal sheaves from $(M\setminus \mathfrak c)/{\mathbb R}_+$.
\end{Remark}

\begin{Remark}
\label{wrapped}
This remark is devoted to a~conjectural parallel picture for the partially wrapped Fukaya category.
Let us continue with the setting of Theorem~\ref{introtheorem2}, and write $WF_\Lambda(M)$ for the partially wrapped Fukaya
category as developed in~\cite{AurICM, Aur}.
It specializes to the fully wrapped variant of~\cite{AS} when the characteristic cone $\Lambda \subset M$ coincides with the
compact core $K\subset M$.

Let $\operatorname{Perf} WF_\Lambda(M)$ denote the stable category of perfect modules over $WF_\Lambda(M)$, and let
$\operatorname{Mod} WF_\Lambda(M)$ denote the stable category of all modules.
Note that the two stabilizations are formally the same amount of information: $\operatorname{Mod} WF_\Lambda(M)$ is the
ind-category $\operatorname {Ind} (\operatorname{Perf} WF_\Lambda(M))$, and $\operatorname{Perf} WF_\Lambda(M)$ is the full
subcategory of compact objects $(\operatorname{Mod} WF_\Lambda(M))^c$.

Given a~conic open subset $U\subset M$, let us imagine a~``mirror'' picture of the category ${\mathcal F}_\Lambda(U)$ as perfect
quasicoherent sheaves over a~scheme $X_U$ proper over $\operatorname{Spec} k$.
The statement that $X_U$ is proper can be formalized by the expectation that ${\mathcal F}_\Lambda(U)$ is hom-f\/inite.
Given conic open subsets $V\subset U\subset M$, let us also imagine a~``mirror'' picture of the restriction map $ {\mathcal
F}_\Lambda(U) \to {\mathcal F}_\Lambda(V) $ as the pullback of perfect quasicoherent sheaves under a~morphism $X_V \to X_U$.

Now given a~conic open subset $U\subset M$, consider the category of f\/inite functionals
\begin{gather*}
\operatorname{Coh} {\mathcal F}_\Lambda(U) = \operatorname {Hom}_{{\rm st}_k}\big({\mathcal F}_\Lambda(U)^{\operatorname{op}},\operatorname{Perf} k\big).
\end{gather*}
From the above ``mirror'' perspective, $\operatorname{Coh} {\mathcal F}_\Lambda(U) $ corresponds to coherent sheaves over the
proper scheme $X_U$.
Given conic open subsets $V\subset U\subset M$, the restriction map $ {\mathcal F}_\Lambda(U) \to {\mathcal F}_\Lambda(V) $
induces a~corestriction map
\begin{gather*}
\xymatrix{ \operatorname{Coh} {\mathcal F}_\Lambda(V) \ar[r] & \operatorname{Coh} {\mathcal F}_\Lambda(U)}.
\end{gather*}
From the above ``mirror'' perspective, the corestriction map corresponds to the pushforward (right adjoint to pullback) of
coherent sheaves under the morphism $X_U \to X_V$.

The above constructions equip $\operatorname{Coh} {\mathcal F}_\Lambda$ with the structure of a~pre-cosheaf.
Let us write $\operatorname {Ind} (\operatorname{Coh} {\mathcal F}_\Lambda)$ for the corresponding pre-cosheaf of
ind-categories, and $\operatorname {Ind}( \operatorname{Coh} {\mathcal F}_\Lambda)^+$ for its cosheaf\/if\/ication.
(If one prefers the intuitions of sheaves over cosheaves, one could pass to the right adjoints of the corestriction maps of
$\operatorname {Ind}( \operatorname{Coh} {\mathcal F}_\Lambda)^+$ and turn it into a~sheaf.) We expect that there is a~natural
equivalence
\begin{gather*}
\xymatrix{ \Gamma(M, \operatorname {Ind}( \operatorname{Coh} {\mathcal F}_\Lambda)^+) \ar[r]^-\sim & \operatorname{Mod}
WF_\Lambda(M) }.
\end{gather*}
In other words, the stable partially wrapped Fukaya category $\operatorname{Perf} WF_\Lambda(M)$ consists of the compact objects
of the global sections of the cosheaf $\operatorname {Ind}( \operatorname{Coh} {\mathcal F}_\Lambda)^+$.

Finally, let us mention some corroborating evidence for the above picture.
First, using the results of~\cite{Abwrapped, N, NZ}, one can check it for $M$ a~cotangent bundle and $\Lambda$ the zero section.
Second, recent calculations of Fukaya categories of Riemann surfaces, pairs-of pants and their genera\-lizations~\cite{AAEKO,Sgenus2, Sheridan} reveal a~mirror symmetry with matrix factorizations.
In the situations considered, the wrapped Fukaya category corresponds to ``coherent'' matrix factorizations and the full
subcategory of compact branes corresponds to ``perfect'' matrix factorizations.
\end{Remark}

We have brief\/ly mentioned anticipated applications of the recollement pattern and Theorem~\ref{introtheorem2} earlier in the
introduction.
Let us conclude here with a~useful technical application: Theorem~\ref{introtheorem2} allows us to def\/ine the stable Fukaya
category of an ``inexhaustible, singular Weinstein manifold'' (for example, an ``open Weinstein cobordism'').
Namely, given a~union of coisotropic cells $C\subset M$, we can take the sections of the sheaf ${\mathcal F}_\Lambda$ lying
strictly above $C$ for increasing $\Lambda$.
Rather than leaving it in such an abstract form, let us explain what results by appealing directly to the concrete geometry of
the recollement pattern.

To begin, let us rotate our viewpoint on the exact triangles appearing in diagram~\eqref{eq triangles} and switch the roles of
the known and unknown categories in our discussion.
Recall that we assumed that the closed coisotropic subvariety $i:C\to M$ is both a~smooth submanifold and a~union of unstable
cells.
Let us relax these two requirements in turn.

First, suppose the closed coisotropic subvariety $i:C\to M$ is no longer necessarily smooth, but still a~union of unstable
cells.
On the one hand, the formalism of Hamiltonian reduction still provides a~correspondence
\begin{gather*}
\xymatrix{ N & \ar@{->>}[l]_-{q} C \ar@{^(->}[r]^-{i} & M }
\end{gather*}
though $N$ is now a~``singular Weinstein manifold'', whatever that might mean.
On the other hand, the open complement $j:M^\circ = M \setminus C\to M$ continues to inherit the structure of a~Weinstein
manifold.
With this setup, the constructions underlying diagram~\eqref{eq functors} should provide adjunctions
\begin{gather*}
\xymatrix{ \mathfrak j^*: \operatorname{Perf} F(M) \ar@<0.5ex>[r] &\ar@<0.5ex>[l] \operatorname{Perf} F(M^\circ):\mathfrak j_*, &
\mathfrak j_!: \operatorname{Perf} F(M^\circ ) \ar@<0.5ex>[r] &\ar@<0.5ex>[l] \operatorname{Perf} F(M):\mathfrak j^! }
\end{gather*}
with $\mathfrak j_!$, $\mathfrak j_*$ fully faithful, and $\mathfrak j^!\simeq \mathfrak j^*$.
We do not have an {\em a~priori} def\/inition of a~Fukaya ca\-te\-gory $\operatorname{Perf} F(N)$, but the recollement pattern tells
us what it should be.
Namely, we should def\/ine~$F(N)$ to classify triples of data~$L \in \operatorname{Perf} F(M)$, $L^\circ \in \operatorname{Perf}
F(M^\circ)$, together with a~morphism $u \in \operatorname {Hom}_{\operatorname{Perf} F(M)}(L, \mathfrak j_*L^\circ)$.
When $N$ is smooth, the recollement pattern conf\/irms that we recover precisely $\operatorname{Perf} F(N)$ by taking the kernel
of morphisms~$u$ appearing in such data.

\begin{Remark}
One should view the preceding as more than a~formal analogue of the situation for ${\mathcal D}$-modules on singular varieties.
There Kashiwara's theorem conf\/irms that such an approach provides an unambiguous notion of ${\mathcal D}$-module.
\end{Remark}

Second, suppose the closed coisotropic subvariety $i:C\to M$ is a~smooth submanifold invariant under the Liouville f\/low, but not
necessarily a~union of unstable cells.
Suppose as well that its Hamiltonian reduction $N$ is a~smooth manifold, and hence a~Weinstein manifold.
Observe that the open complement $j:M^\circ = M \setminus C\to M$ may be viewed as an ``inexhaustible Weinstein manifold'', or
informally speaking, a~complete exact symplectic manifold such that the Liouville vector f\/ield is gradient-like for a~possibly
inexhausting Morse function.
With this setup, the constructions underlying diagram~\eqref{eq functors} should provide adjunctions
\begin{gather*}
\xymatrix{ \mathfrak i^*: \operatorname{Perf} F(N) \ar@<0.5ex>[r] &\ar@<0.5ex>[l] \operatorname{Perf} F(M):\mathfrak i_*, &
\mathfrak i_!: \operatorname{Perf} F(M) \ar@<0.5ex>[r] &\ar@<0.5ex>[l] \operatorname{Perf} F(N):\mathfrak i^! }
\end{gather*}
with $\mathfrak i_!\simeq\mathfrak i_*$ fully faithful.
We do not have an {\em a~priori} def\/inition of a~Fukaya category $\operatorname{Perf} F(M^\circ)$, but the recollement pattern
tells us what it should be.
Namely, we should def\/ine $\operatorname{Perf} F(M^\circ)$ to classify triples of data $L_N \in \operatorname{Perf} F(N)$, $L \in
\operatorname{Perf} F(M)$, together with a~morphism $c \in \operatorname {Hom}_{\operatorname{Perf} F(M)}(\mathfrak i_! L_N,
L)$.
When $C$ is a~union of unstable cells, the the recollement pattern conf\/irms that we recover precisely $\operatorname{Perf}
F(M^\circ)$ by taking the cokernel of morphisms $c$ appearing in such data.

Putting together the above generalizations, we obtain an unambigious inf\/initesimal Fukaya category of an ``inexhaustible,
singular Weinstein manifold'' so that it is compatible with familiar notions whenever they apply.

\subsection{Inf\/luences}

In the remainder of the introduction, we recount some inf\/luences on our thinking, in particular the symplectic geometry of Dehn
twists, the recollement pattern for constructible sheaves, and the Morse theory of integral kernels.

\subsubsection{Dehn twists}

One can view our expectations as simple elaborations on the fundamental notion of exact symplectic Dehn twists.
In turn, the many incarnations of this notion (spherical twists, mutations, braid actions, Hecke correspondences, \dots)
play a~prominent role in mirror symmetry and geometric representation theory.

We sketch here an informal picture of our expectations from the perspective of Dehn twists.
The basic mantra could be: {\em Dehn twists around spheres provide mutations; Dehn twists around cells provide semiorthogonal
decompositions}.
(There are also extensive relations between mutations and semiorthogonal decompositions, often arising from exceptional
collections, thanks to the fact that spheres themselves can be cut into cells.)

Let us recall Seidel's long exact sequence in Floer cohomology~\cite{S01b,S01, S03}.
There are also highly relevant relative sequences found in the work of Perutz~\cite{P1, P2, P08} and
Wehrheim--Woodward~\cite{ww6}.
We will ignore technical issues and proceed as quickly as possible to the statement.

Let $S\subset M$ be an exact Lagrangian sphere in an exact symplectic manifold.
Let $\tau_S:M\to M$ be the associated exact symplectic Dehn twist around $S$.
Then for any two exact Lagrangian submanifolds $L_0, L_1 \subset M$, there is a~long exact sequence of (${\mathbf Z}/2{\mathbf
Z}$-graded, ${\mathbf Z}/2{\mathbf Z}$-linear) Floer cohomology groups
\begin{gather*}
\xymatrix{ HF(\tau_S(L_0), L_1) \ar[r] & HF(L_0, L_1) \ar[r] & HF(S, L_1) \otimes HF(L_0, S) \ar[r]^-{[1]} & }.
\end{gather*}

The sequence admits a~straightforward categorical interpretation.
Assume $M$ is equipped with appropriate background structures, and the exact Lagrangian submanifolds $S$, $L_0$, $L_1$ are all
equipped with appropriate brane structures.
Let $\operatorname{Mod} F(M)$ denote the ${\mathbf Z}$-graded, $k$-linear stable Fukaya category of modules, and let $S$, $L_0$,
$L_1$ denote the corresponding objects.
Thanks to the functoriality of the sequence in the variable $L_1$, we can rewrite it as an exact triangle
\begin{gather*}
\xymatrix{ S \otimes \operatorname {Hom}_{\operatorname{Mod} F(M)}(S, L_0) \ar[r] & L_0 \ar[r] & \tau_S(L_0) \ar[r]^-{[1]} & }.
\end{gather*}

Let us introduce the (presently elaborate but ultimately justif\/ied) notation
\begin{gather*}
\xymatrix{ \mathfrak i_{S!}:\operatorname{Mod} F_S(M) \ar@<0.5ex>[r] &\ar@<0.5ex>[l] \operatorname{Mod} F(M):\mathfrak i^!_S }
\end{gather*}
for the fully faithful embedding $\mathfrak i_{S!}$ of the subcategory $\operatorname{Mod} F_S(M)$ generated by $S$, and its
right adjoint $\mathfrak i^!_S = \operatorname {Hom}_{\operatorname{Mod} F(M)}(S, -)$.
Then thanks to the functoriality of the sequence in the va\-riab\-le~$L_0$, we can view it as an exact triangle of functors
\begin{gather*}
\xymatrix{ \mathfrak i_{S!} \mathfrak i_S^! \ar[r] & \operatorname{id}_{\operatorname{Mod} F(M)} \ar[r] & \tau_S \ar[r]^-{[1]}& },
\end{gather*}
where the f\/irst map is the counit of the adjunction.

Now suppose the exact Lagrangian sphere $S\subset M$ were rather a~closed but noncompact exact Lagrangian cell $C\subset M$.
To place it in a~categorical context, let us now allow the Fukaya category $F(M)$ to contain closed but noncompact Lagrangian
branes.
Because $C \subset M$ is now a~cell, we will be able to work with the ${\mathbf Z}$-graded, $k$-linear stable category
$\operatorname{Perf} F(M)$ of perfect modules.
Then one expects to have a~similar exact triangle of functors
\begin{gather*}
\xymatrix{ \mathfrak i_{C!} i_C^! \ar[r] & \operatorname{id}_{\operatorname{Perf} F(M)} \ar[r] & \tau_C \ar[r]^-{[1]} & },
\end{gather*}
where $\tau_C$ denotes the ``Dehn twist'' around the cell $C\subset M$.
It takes compact Lagrangian branes to noncompact Lagrangian branes which are asymptotically close to $C$ near inf\/inity.

Let us go one step further and consider the open exact symplectic submanifold $\mathfrak j_{M^\circ}:M^\circ = M \setminus C \to
M$.
Assuming $C$ is in good position with respect to the exact symplectic structure, we should be able to realize $\tau_C$ as the
monad of an adjunction
\begin{gather*}
\xymatrix{ \mathfrak j_{M^\circ}^*: \operatorname{Perf} F(M) \ar@<0.5ex>[r] &\ar@<0.5ex>[l] \operatorname{Perf} F(M^\circ
):\mathfrak j_{M^\circ*} },
\end{gather*}
where $j_{M^\circ*}$ is a~fully faithful embedding.
Here the functor $j_{M^\circ*}$ is geometric, built out of the inclusion $j_{M^\circ}$ and the Dehn twist $\tau_C$ near
inf\/inity.
Its left adjoint $j_{M^\circ}^*$ is of a~categorical origin, just as the left adjoint $i^!_C$ is of a~categorical origin.

Putting the above together, we obtain an exact triangle of functors
\begin{gather*}
\xymatrix{ \mathfrak i_{C!} \mathfrak i_C^! \ar[r] & \operatorname{id}_{\operatorname{Perf} F(M)} \ar[r] &\mathfrak
j_{M^\circ*}\mathfrak j_{M^\circ}^* \ar[r]^-{[1]} & },
\end{gather*}
where the initial map is the counit of the adjunction, and the middle map is the unit of the adjunction.
Thus we have a~semiorthogonal decomposition of $\operatorname{Perf} F(M)$ by the two full subcategories $\operatorname{Perf}
F_C(M)$, $\operatorname{Perf} F(M^\circ)$ with the semiorthogonality $\operatorname{Perf} F_C(M)^\perp \simeq \operatorname{Perf}
F(M^\circ)$, $\operatorname{Perf} F_C(M) \simeq {}^\perp \operatorname{Perf} F(M^\circ)$.

We will suggest a~generalization of the above picture where we allow the cell $C\subset M$ to be coisotropic rather than
Lagrangian.
The geometric constructions and formal consequences should be similar, with the main new development that the full subcategory
$\operatorname{Perf} F_C(M)$ no longer should be generated by a~single object.
In the section
immediately following, we motivate the general pattern with the formalism of recollement for constructible sheaves.

\subsubsection{Recollement pattern}
\label{dev}

We recall here the recollement pattern for constructible sheaves.

In what follows, we will only consider subanalytic sets $X$ and subanalytic maps $f:X\to Y$.
We write $\operatorname{Sh}(X)$ for the dif\/ferential graded category of constructible complexes of sheaves on $X$, and have the
standard adjunctions
\begin{gather*}
\xymatrix{ f^*:\operatorname{Sh}(Y) \ar@<0.5ex>[r] &\ar@<0.5ex>[l] \operatorname{Sh}(X):f_*,
& f_!:\operatorname{Sh}(X)\ar@<0.5ex>[r] &\ar@<0.5ex>[l] \operatorname{Sh}(Y):f^! }.
\end{gather*}
Verdier duality provides an anti-involution
\begin{gather*}
\xymatrix{ {\mathbb D}_X:\operatorname{Sh}(X) \ar[r]^-\sim & \operatorname{Sh}(X)^{\operatorname{op}} }
\end{gather*}
intertwining the preceding adjunctions
\begin{gather*}
f^! \simeq {\mathbb D}_X f^* {\mathbb D}_Y,
\qquad
f_! = {\mathbb D}_Y f_* {\mathbb D}_X.
\end{gather*}

Suppose we have a~partition of $X$ into an open subset and its closed complement
\begin{gather*}
\xymatrix{ j:U\ar@{^(->}[r] & X & \ar@{_(->}[l] Y = X\setminus U: i }.
\end{gather*}
Then the standard functors provide a~diagram
\begin{gather*}
\xymatrix{ \operatorname{Sh}(U) \ar@/^2pc/[rr]^-{j_!}\ar@/_2pc/[rr]^-{j_*} && \ar[ll]_-{j^! \simeq j^*} \operatorname{Sh}(X)
\ar@/^2pc/[rr]^-{i^*} \ar@/_2pc/[rr]^-{i^!} && \ar[ll]_-{i_! \simeq i_*} \operatorname{Sh}(Y) }
\end{gather*}

The fact that $j$ is an open embedding (hence smooth of relative dimension zero) provides a~canonical identif\/ication $j^! \simeq
j^*$, and the fact that $i$ is a~closed embedding (hence proper) provides a~canonical identif\/ication $i_! \simeq i_*$.

The fact that $U$ and $Y$ are disjoint provide canonical identif\/ications $j^*\mathfrak i_! \simeq 0 \simeq\mathfrak j^! i_*$ and
$i^*j_! \simeq 0 \simeq i^! j_* $.
The additional fact that $U$ and $Y$ cover $X$ lead to dual exact triangles
\begin{gather*}
\xymatrix{i_{!} i^! \ar[r] & \operatorname{id}_{\operatorname{Sh}(X)} \ar[r] & j_*j^* \ar[r]^-{[1]} &, &  j_{!}\mathfrak j^!
\ar[r] & \operatorname{id}_{\operatorname{Sh}(X)} \ar[r] & i_*i^* \ar[r]^-{[1]} & } \!\!,
\end{gather*}
where the f\/irst and middle morphisms are respectively counits and units of adjunctions.
These triangles are generalizations of the long exact sequences of pairs in cohomology.

Alternatively, we also have the dual exact triangles
\begin{gather*}
\xymatrix{\!\!\operatorname{id}_{\operatorname{Sh}(X)}\ar[r] & j_* j^* \oplus i_*i^* \ar[r] & i_*i^* j_* j^* \ar[r]^-{[1]}& \!,  &
\mathfrak i_!i^!\mathfrak j_!\mathfrak j^! \ar[r] & \mathfrak j_!\mathfrak j^! \oplus\mathfrak i_!i^! \ar[r] &
\operatorname{id}_{\operatorname{Sh}(X)}\ar[r]^-{[1]} &}\!\! ,
\end{gather*}
where the f\/irst and middle morphisms result from units and counits of adjunctions respectively.
These triangles are generalizations of the Mayer--Vietoris long exact sequences of cohomology.

Now let us f\/ind the above formalism as a~special case of the recollement pattern of diagrams~\eqref{eq functors} and~\eqref{eq
triangles}.
By inductive considerations, one can see that the above formalism devolves from the case when $X$ is a~manifold and $Y\subset X$
is a~submanifold.
Let us focus on this case and consider the Weinstein manifolds $M = T^*X$ and $N = T^*Y$.
Recall from~\cite{N, NSpr, NZ} the microlocalization equivalence
\begin{gather*}
\xymatrix{ \mu_X:\operatorname{Sh}(X) \ar[r]^-\sim & \operatorname{Perf} F(T^*X) }.
\end{gather*}

Then by an inductive sequence of applications of the recollement pattern of of diagrams~\eqref{eq functors} and~\eqref{eq
triangles}, we obtain commutative diagrams with horizontal maps fully faithful embeddings and vertical maps equivalences
\begin{gather*}
\xymatrix{ \ar[d]_-{\mu_Y}\operatorname{Sh}(Y) \ar@{^(->}[r]^{i_!} & \operatorname{Sh}(X)\ar[d]^-{\mu_X} & \ar[d]_-{\mu_U}
\operatorname{Sh}(U ) \ar@{^(->}[r]^-{j_*} & \operatorname{Sh}(X)\ar[d]^-{\mu_X}
\\
\operatorname{Perf} F(T^*Y) \ar@{^(->}[r]^{\mathfrak i_!} & \operatorname{Perf} F(T^*X) & \operatorname{Perf} F(T^*U )
\ar@{^(->}[r]^-{\mathfrak j_*} & \operatorname{Perf} F(T^*X) }
\end{gather*}

Using the natural involutive anti-symplectomorphism on a~cotangent bundle, one can construct Verdier duality on its Fukaya
category as well~\cite{N}.

\subsubsection{Morse theory of integral kernels}

Our expectations are guided by the well-known strategy: to prove universal statements about objects of a~category, one should
realize endofunctors of the category as integral transforms and establish canonical identities among them.

We will illustrate this with the toy case of Morse theory since our later discussion is a~direct analogue of it.
For simplicity, let us work with a~compact oriented manifold $M$.
Consider a~generic pair consisting of a~Morse function $f:M\to {\mathbb R}$ and Riemannian metric $g$ on $M$.
Let $\mathfrak c\subset M$ denote the critical locus of $f$, and $\Phi_t:M \to M$ the f\/low along the gradient $\nabla_g f$.
To each critical point $p\in\mathfrak c$, associate the stable and unstable cells
\begin{gather*}
S_p = \big\{ x\in M \, |\, \lim_{t\to \infty} \Phi_t(x) = p\big\},
\qquad
U_p = \big\{ x\in M \, |\, \lim_{t\to -\infty} \Phi_t(x)=p\big\}.
\end{gather*}
Finally, orient $S_p$ and $U_p$ so that the orientation of $U_p\times S_p$ agrees with that of $X$ at $p$.

One formulation of Morse theory is that every cohomology class $c\in H^*(M; {\mathbb R})$ can be expressed in the form
\begin{gather}
\label{morse}
c = \sum\limits_{p\in \mathfrak c} \langle c, U_p\rangle S_p,
\end{gather}
where $\langle c_p, U_p\rangle $ denotes the natural pairing, and we regard $S_p$ as a~cohomology class via Poincar\'e duality.
In other words, the stable and unstable cells form dual bases in cohomology.

To establish equation~\eqref{morse}, it is useful to recast it as an equation in the cohomology of the product $X\times X$.
Namely, each cohomology class $ k\in H^*(X\times X; {\mathbb R})$ can be regarded as an integral kernel providing an
endomorphism
\begin{gather*}
\xymatrix{ \Phi_k:H^*(X;{\mathbb R})\ar[r] & H^*(X; {\mathbb R}), & \Phi_k(c) = p_{2!}(p_1^*(c) \cap k) },
\end{gather*}
where we use Poincar\'e duality to integrate.
Then equation~\eqref{morse} follows from the identity of cohomology classes
\begin{gather}
\label{morsekernels}
\Delta_X = \sum\limits_{p\in \mathfrak c} U_p\times S_p,
\end{gather}
where $\Delta_X \subset X \times X$ is the diagonal, and hence provides the identity endomorphism.
Finally, to establish equation~\eqref{morsekernels}, one observes that the gradient f\/low of the Morse function
\begin{gather*}
\xymatrix{ f\circ p_1 - f\circ p_2:X\times X\ar[r] & {\mathbb R} }
\end{gather*}
provides a~homotopy between the diagonal $\Delta_X$ and the sum of external products $\sum\limits_{p\in \mathfrak c} U_p\times
S_p$.

All in all, a~pleasant aspect of the above argument is that it makes no reference to an arbitrary cohomology class $c\in H^*(X;
{\mathbb R})$, but only involves highly structured cohomology classes such as the diagonal $\Delta_X$ and the sum of external
products $\sum\limits_{p\in \mathfrak c} U_p\times S_p$.
Furthermore, we do not need to know a~precise relation between cohomology classes $k \in H^*(X\times X; {\mathbb R})$ and
arbitrary endomorphisms $\Phi:H^*(X;{\mathbb R}) \to H^*(X;{\mathbb R})$, only that there is a~linear map
\begin{gather*}
\xymatrix{ H^*(X\times X; {\mathbb R})\ar[r] & \operatorname{End} H^*(X; {\mathbb R}), & k \ar@{|->}[r] & \Phi_k }.
\end{gather*}

Turning from vector spaces to the setting of linear categories, one f\/inds many examples of the above argument.
Most prominently, there is Beilinson's universal resolution~\cite{Beil} of coherent sheaves on projective space $\mathbb P^n$ by
vector bundles.
It suf\/f\/ices once and for all to introduce the Koszul resolution of the structure sheaf ${\mathcal O}_{\Delta_{\mathbb P^n}}$ of
the diagonal $\Delta_{\mathbb P^n} \subset \mathbb P^n\times \mathbb P^n$.
Then for any coherent sheaf on $\mathbb P^n$, convolution with the Koszul resolution produces the desired resolution by vector
bundles.
What results is a~concrete description of all coherent sheaves in terms of the quiver of constituent vector bundles in the
Koszul resolution.

Our sketched evidence for the recollement pattern of diagrams~\eqref{eq functors} and~\eqref{eq triangles} applies the above
version of Morse theory in the setting of linear categories.
We consider the product Weinstein manifold $M^{\operatorname{op}} \times M$ where we write $M^{\operatorname{op}}$ to denote the opposite symplectic
structure.
The formalism of bimodules should allow us to view Lagrangian branes $L \subset M^{\operatorname{op}}\times M$ as endofunctors.
In particular, the diagonal brane $\Delta_M \subset M^{\operatorname{op}} \times M$ should represent the identity functor.
The endofunctors of adjunctions appearing in the exact triangles of~\eqref{eq triangles} should also be represented by natural
correspondences.

A familiar but key observation is that the geometry of a~neighborhood of the diagonal $\Delta_M$ looks like the geometry of the
cotangent bundle $T^*M$.
We highlight a~useful elaboration: the geometry of natural correspondences should match the geometry of the conormal bundles to
the coisotropic unstable cells.

\section{Weinstein manifolds}
\label{sect weinstein}

Most of the material of this section is well-known and available in many beautiful sources~\cite{CE, E, EG, W}.
Unless otherwise stated, we will assume that all manifolds are real analytic, and all subsets and maps are
subanalytic~\cite{BieMi, vdDM}.

\subsection{Basic notions}

\begin{Definition}
An {\em exact symplectic manifold} $(M, \theta)$ is a~manifold $M$ with a~one-form $\theta$ such that $\omega = d\theta$ is
a~symplectic form.
We use the term {\em Liouville form} to refer to $\theta$.
\end{Definition}

One def\/ines the Liouville vector f\/ield $Z$ by the formula $i_Z \omega = \theta$.
The symplectic form $\omega$ is an eigenvector for the Lie derivative $L_Z \omega = i_Z d\omega + d i_Z\omega = 0 + d\theta =
\omega$.
An exact symplectic manifold is equivalently a~triple $(M, \omega, Z)$ consisting of a~manifold $M$ with symplectic form
$\omega$ and f\/ixed vector f\/ield $Z$ such that $L_Z\omega = \omega$.

We will always assume that $M$ is {\em complete} with respect to $Z$ in the sense that $Z$ integrates for all time to provide an
expanding action
\begin{gather*}
\xymatrix{ \Phi_{t}:{\mathbb R}_{+ } \times M \ar[r] & M }.
\end{gather*}

\begin{Definition}
A subset $A\subset M$ is said to be {\em conic} if it is invariant under the expanding action in the sense that $\Phi_{t}(A) =
A$, for all $t\in {\mathbb R}_+$.
\end{Definition}

\begin{Example}
Any manifold $X$ provides an exact symplectic manifold $(T^*X, \theta_X)$ given by the cotangent bundle $\pi_X:T^*X\to X$
equipped with its canonical exact structure $\theta_X$.
The Liouville vector f\/ield $Z_X$ generates the standard linear scaling $\Phi_{X, t}$ along the f\/ibers of $\pi_X$.

More generally, any conic open subset $A\subset T^*X$ in particular the complement of the zero section $T^*X\setminus X$,
provides an exact symplectic manifold $(A, \theta_X|_A)$.
\end{Example}

By a~Morse function $h:M \to {\mathbb R}$, we will always mean an exhausting (proper and bounded below) function whose critical
points are nondegenerate and f\/inite in number.
A vector f\/ield~$V$ is gradient-like with respect to $h$ if away from the critical points of~$h$, we have $dh(V)>0$, and in some
neighborhood of the critical points, $V$ is the gradient of $h$ with respect to some Riemannian metric.

\begin{Definition}
By a~{\em Weinstein manifold}, we will mean an exact symplectic manifold $(M, \theta)$ that admits a~Morse function $h:M\to
{\mathbb R}$ such that the Liouville vector f\/ield $Z$ is gradient-like with respect to $h$.
\end{Definition}

\begin{Example}
A {\em Weinstein cell} is a~Weinstein manifold $(M, \theta)$ such that $\theta$ has a~single zero.
It follows that any Morse function $h:M\to {\mathbb R}$ such that the Liouville vector f\/ield $Z$ is gradient-like with respect
to $h$ will have a~single critical point which is a~minimum.
\end{Example}

\begin{Example}
\label{cotangent is weinstein}
Suppose $X$ is a~compact manifold.
Choose a~generic Morse function \mbox{$f_X:X\to {\mathbb R}$} and Riemannian metric $g$, and let $\nabla_g f_X$ denote the resulting
gradient of~$f_X$.
We obtain a~f\/iberwise linear function $F_X:T^*X \to {\mathbb R}$ by setting $F_X(x, \xi) = \xi(\nabla_g f_X|_x)$.

Then for $\epsilon>0$ suf\/f\/iciently small, the pair $(T^*X, \theta_X + \epsilon dF_X)$ forms a~Weinstein manifold, exhibited by
the Morse function $h = g + \pi^*_Xf_X: T^*X\to {\mathbb R}$ given by $h(x, \xi) = |\xi|_g^2 + f_X(x)$.
\end{Example}

\subsection{Cell decompositions} For the statements of this section, in addition to the previously mentioned works, we recommend
the excellent sources~\cite{biran01, biran06}.

Let $(M, \theta)$ be a~complete exact symplectic manifold with Liouville vector f\/ield $Z$ and expan\-ding action
\begin{gather*}
\xymatrix{ \Phi_t:{\mathbb R}_{+ } \times M \ar[r] & M }.
\end{gather*}

\begin{Definition}
The {\em critical locus } $\mathfrak c\subset M$ is the zero-locus of the exact structure
\begin{gather*}
\xymatrix{ \mathfrak c = \{ x\in M \, |\, \theta|_{T_x M} = 0\}.
}
\end{gather*}

The {\em core} $K\subset M$ is the subset
\begin{gather*}
  K = \big\{ x\in M \, |\, \lim_{t\to \infty} \Phi_t(x) \in \mathfrak c \big\}.
\end{gather*}

The {\em ether} $E\subset M$ is the complement
\begin{gather*}
E = M\setminus K.
\end{gather*}
\end{Definition}

\begin{Remark}
All three of the above subsets are evidently conic, and the expanding action on the complement $M\setminus \mathfrak c$, and in
particular the ether $E$, is free.

The critical locus $\mathfrak c\subset M$ is equivalently the zero-locus of the Liouville vector f\/ield $Z$, or the f\/ixed points
of the expanding action $\Phi_t$.
\end{Remark}

\begin{Example}
For the exact symplectic manifold $(T^*B, \theta_B)$, the core and critical locus are the zero section $\mathfrak c_{B} = K_{B}
= B $, and the ether is its complement $E_{B} = T^*B \setminus B$.

More generally, for any conic open subset $A\subset T^*B$ and resulting exact symplectic manifold $(A, \theta_B|_A)$, its
critical locus and core are the intersection $\mathfrak c = K = A\cap B$, and its ether is the complement $E = A \cap (T^*B
\setminus B)$.
\end{Example}

Now suppose $(M, \theta)$ is a~Weinstein manifold.
Then the critical locus $\mathfrak c \subset M$ consists of the f\/initely many critical points of a~Morse function $h:M \to
{\mathbb R}$ for which $Z$ is gradient-like.

\begin{Definition}
Given a~Weinstein manifold $(M, \theta)$, to each critical point $p\in\mathfrak c$, we associate the {\em stable} and {\em
unstable manifolds}
\begin{gather*}
  c_p = \big\{ x\in M \, |\, \lim_{t\to \infty} \Phi_t(x) = p\big\}, \qquad C_p = \big\{ x\in M \, |\, \lim_{t\to -\infty} \Phi_t(x)=p\big\} ,
\end{gather*}
\end{Definition}

\begin{Lemma}
Given a~Weinstein manifold $(M, \theta)$, for each critical point $p\in\mathfrak c$, the stable manifold $c_p$ is an isotropic
cell and the unstable manifold $C_p$ is a~coisotropic cell.
\end{Lemma}

The coisotropic cells provide a~partition
\begin{gather*}
M = \coprod_{p\in \mathfrak c} C_p
\end{gather*}
while the isotropic cells partition the core
\begin{gather*}
K = \coprod_{p\in \mathfrak c} c_p \subset M.
\end{gather*}

We will always perturb our Liouville form $\theta$ to be generic.
Then it follows that the partition by coisotropic cells is a~Whitney stratif\/ication.

\begin{Lemma}
Given a~Weinstein manifold $(M, \theta)$, the core $K\subset M$ is compact and its inclusion is a~homotopy-equivalence.
\end{Lemma}

We will use the coisotropic cells to view the f\/inite critical locus $\mathfrak c \subset M$ as a~{poset}: we say critical points
$p, q\in \mathfrak c$ satisfy $p \leq q$ if and only if
\begin{gather*}
\overline C_p \cap C_q \not = \varnothing.
\end{gather*}
In particular, maxima $p\in \mathfrak c$ correspond to closed coisotropic cells $C_p\subset M$.
Subsets $\mathfrak s\subset \mathfrak c$ such that $p\in \mathfrak s$ and $q\in \mathfrak c$, with $q\leq p$, implies $q\in
\mathfrak s$ correspond to open unions
\begin{gather*}
M_\mathfrak s = \coprod_{p\in \mathfrak s} C_p \subset M.
\end{gather*}

\subsection{Markings}

Let us f\/irst list some useful notions about conic sets in a~complete exact symplectic manifold $(M, \theta)$.

Note that the Liouville f\/low $\Phi_t:{\mathbb R}_+\times M\to M$ is free on the ether $E\subset M$ and so the quotient
$E/{\mathbb R}_+$ is a~manifold.
Furthermore, since $\theta$ vanishes on the Liouville vector f\/ield $Z$, and the Lie derivative satisf\/ies $L_Z\theta = \theta$,
it makes sense to ask whether $\theta$ is positive, zero or negative on a~tangent vector to $E/{\mathbb R}_+$.
It is easily checked that $\ker(\theta) \subset T(E/{\mathbb R}_+)$ provides a~canonical co-oriented contact structure.

\begin{Definition}
Let $(M, \theta)$ be a~complete exact symplectic manifold.

(1) The {\em projectivization} of $M$ is the contact manifold
\begin{gather*}
M^\infty = E/{\mathbb R}_{+}
\end{gather*}
equipped with its canonical contact structure $\xi = \ker(\theta)$.

More generally, the {\em projectivization} of a~conic subset $A\subset M$ is the subset
\begin{gather*}
A^\infty = (A\cap E)/{\mathbb R}_{+} \subset M^\infty.
\end{gather*}

(2) The {\em cone} over a~subset $A^\infty \subset M^\infty$ is the conic subset
\begin{gather*}
\operatorname{c}A^\infty =\big\{x\in E \, |\, \lim_{t\to \infty} \Phi_t(x) \in A^\infty\big\}  \subset M.
\end{gather*}
\end{Definition}

\begin{Remark}
When $(M, \theta)$ is a~Weinstein manifold, its projectivization $M^\infty$ is compact and contactomorphic to any level set
$h^{-1}(r)\subset M$ of a~compatible Morse function $h:M\to \mathbb R$, for $r$ greater than all critical values.
\end{Remark}

Next let us list some useful notions about closed sets.

\begin{Definition}
Let $(M, \theta)$ be a~complete exact symplectic manifold.

(1) The {\em compactification at infinity} of $M$ is the manifold with boundary
\begin{gather*}
\overline M = \left((M \times (0, \infty])\setminus (K \times \{\infty\})\right)/ {\mathbb R}_{+} \simeq M \cup M^\infty.
\end{gather*}

More generally, the {\em compactification at infinity} of a~closed subset $A\subset M$ is the closure
\begin{gather*}
\overline A \subset \overline M.
\end{gather*}

(2) The {\em boundary at infinity} of a~closed subset $A \subset M$ is the frontier
\begin{gather*}
\partial^\infty A = \overline A \setminus A \subset M^\infty = \overline M \setminus M.
\end{gather*}
\end{Definition}

\begin{Remark}
If $A\subset M$ is closed and conic, then our notation agrees in that $A^\infty = \partial^\infty A$ as subsets of $M^\infty$.
\end{Remark}

Finally, we introduce the notion of a~characteristic cone $\Lambda\subset M$, which could alternatively be called a~``noncompact
skeleton''.
Along with satisfying natural properties stated in the def\/inition below, it will be a~not necessarily smooth subvariety.
By a~subvariety $Y \subset M$, we mean nothing more than a~subanalytic subset of $M$ (or if more generality is desired,
a~def\/inable subset within some o-minimal context).
We will only be interested in $Y$ as a~topological subspace of $M$, rather than with any of its potential algebraic aspects.
In particular, $Y$ admits a~Whitney stratif\/ication by submanifolds of $M$, and following standard convention, we say $Y$ is
isotropic if the symplectic form vanishes when restricted to each stratum.

\begin{Definition}
By a~{\em marked exact symplectic manifold} $(M, \theta, \Lambda)$, we will mean an exact symplectic manifold $(M, \theta)$
together with a~closed conic isotropic subvariety $\Lambda \subset M$ containing the core $K \subset \Lambda$.
We use the term {\em characteristic cone} to refer to $\Lambda$.

By a~{\em marked Weinstein manifold} $(M, \theta, \Lambda)$, we will mean a~marked exact symplectic manifold such that the
underlying exact symplectic manifold $(M, \theta)$ is a~Weinstein manifold.
\end{Definition}

\begin{Remark}
Given a~marked exact symplectic manifold $(M, \theta, \Lambda)$, we can alternatively encode the characteristic cone
$\Lambda\subset M$ by taking its projectivization $\Lambda^\infty\subset M^\infty$.
We recover the characteristic cone $\Lambda\subset M$ by taking the union of the core and the cone over the projectivization
\begin{gather*}
\Lambda = K \cup c\Lambda^\infty \subset M.
\end{gather*}
\end{Remark}

\subsection{Coisotropic cells}

Suppose $(M, \theta, \Lambda)$ is a~marked Weinstein manifold.
Fix a~critical point $p\in \mathfrak c$, and consider the inclusion of the coisotropic cell $i_p:C_p \to M$.

The linear geometry of the inclusion of the coisotropic cell $i_p:C_p \to M$ provides a~useful guide to keep in mind.

First, we have the normal bundle $\mathfrak N i_p \to C_p$ appearing in the exact sequence
\begin{gather*}
\xymatrix{ 0 \ar[r] & TC_p \ar[r] & i_p^*TM \ar[r] & \mathfrak N i_p \ar[r] & 0 }.
\end{gather*}
Dually, we have the conormal bundle $\mathfrak N^* i_p =T^*_{C_p} M\to C_p$ appearing in the exact sequence
\begin{gather*}
\xymatrix{ 0 & \ar[l] T^*C_p & \ar[l] i_p^*T^*M & \ar[l] \mathfrak N^* i_p & \ar[l] 0 }.
\end{gather*}

The symplectic form $\omega = d\theta$ provides an integrable isotropic foliation $\mathfrak f_p\subset TC_p$ by either taking
the symplectic orthogonal $\mathfrak f_p = (TC_p)^\perp$ or equivalently, the symplectic partner $i_{\mathfrak f_p}\omega =
\mathfrak N^* i_p$.
Thus the symplectic form identif\/ies the partial f\/lags
\begin{gather*}
\xymatrix{ \mathfrak f_p \ar@{^(->}[r] & TC_p \ar@{^(->}[r] & i_p^* TM, & \mathfrak N^* i_p \ar@{^(->}[r] & T^*_{\mathfrak f_p} M
\ar@{^(->}[r] & i_p^*T^*M },
\end{gather*}
where $T^*_{\mathfrak f_p} M \subset i_p^* T^*M$ denotes the subbundle annihilating $\mathfrak f_p \subset TC_p \subset i^*_p TM$.
Dually, it identif\/ies the quotient sequences
\begin{gather*}
\xymatrix{ \mathfrak f_p^* & \ar@{->>}[l] T^*C_p & \ar@{->>}[l] i_p^*T^*M, & \mathfrak N i_p & \ar@{->>}[l] i_p^*TM /\mathfrak
f_p & \ar@{->>}[l] i_p^*TM }.
\end{gather*}

Now consider the Hamiltonian reduction diagram
\begin{gather*}
\xymatrix{ M_p & \ar@{->>}[l]_-{q_p} C_p \ar@{^(->}[r]^-{i_p} & M },
\end{gather*}
where $i_p$ is the inclusion of the coisotropic cell, and $q_p$ is the quotient by the integrable isotropic foliation $\mathfrak
f_p \subset TC_p$ determined by $i_p$.
We use the term Hamiltonian reduction diagram with the following constructions in mind.

To f\/ix notation, suppose $\dim C_p = n+k$ where $\dim M= 2n$.
Then we can choose $(n-k)$ local coordinates in a~neighborhood of $p$ such that the intersection of $ C_p$ with the neighborhood
is their zero-locus.
Using the Liouville f\/low $\Phi_t$ we can extend the local coordinates to $(n-k)$ functions in a~neighborhood of $C_p$ such that
$C_p$ is their zero-locus.
These $(n-k)$ functions provide a~moment map for an action of ${\mathbb R}^k$ such that $M_p$ is the quotient symplectic
manifold.
To conf\/irm that the leaf space is indeed a~manifold, and hence a~symplectic manifold, it suf\/f\/ices to check that the foliation is
f\/ibered (or admits local slices).
Observe that this is true locally in $C_p$, in particular in a~neighborhood of $p$.
Then observe that the leaves of the foliation are preserved by the Liouville f\/low $\Phi_t$, and the f\/low of any neighborhood of
$p$ exhausts all of $C_p$.

With the preceding in hand, the rest of the following assertion is evident by construction.

\begin{Lemma}
We have a~natural marked Weinstein cell $(M_p, \theta_p, \Lambda_p)$ characterized by
\begin{gather*}
q_p^*\theta_p = \theta|_{C_p},
\qquad
\Lambda_p = q_p(\Lambda \cap C_p).
\end{gather*}
\end{Lemma}

\begin{Remark}
In the case when $\Lambda^\infty = \varnothing$ so that $\Lambda = K$, we have $\Lambda_p = q_p(K \cap C_p)$.
If in addition $p\in \mathfrak c$ is maximal, so that $C_p\subset M$ is closed, we have $K \cap C_p = \{p\}$, and hence
$\Lambda_p = \{q_p(p)\}$, and so $\Lambda_p^\infty = \varnothing$.
\end{Remark}

\subsection{Recollement cone}

We continue with $(M, \theta, \Lambda)$ a~marked Weinstein manifold.

Fix a~critical point $p\in \mathfrak c$, and return to the Hamiltonian reduction diagram
\begin{gather*}
\xymatrix{ M_p & \ar@{->>}[l]_-{q_p} C_p \ar@{^(->}[r]^-{i_p} & M },
\end{gather*}
where $i_p$ is the inclusion of the coisotropic cell, and $q_p$ is the quotient by the integrable isotropic foliation $\mathfrak
f_p \subset TC_p$ determined by $i_p$.

Recall the natural marked Weinstein cell $(M_p, \theta_p, \Lambda_p)$ characterized by
\begin{gather*}
q_p^*\theta_p = \theta|_{C_p},
\qquad
\Lambda_p = q_p(\Lambda \cap C_p).
\end{gather*}

Let us observe that the inverse-image
\begin{gather*}
\Lambda_{\tilde p} = q_p^{-1}(\Lambda_p) \subset C_p
\end{gather*}
is a~conic isotropic subvariety such that
\begin{gather*}
\Lambda\cap C_p \subset \Lambda_{\tilde p}.
\end{gather*}

\begin{Definition}
The {\em local recollement cone} $\Lambda_{p+} \subset M$ is the conic isotropic subvariety
\begin{gather*}
\Lambda_{p+} = \Lambda \cup \Lambda_{\tilde p} \subset M.
\end{gather*}

The {\em global recollement cone} $\Lambda_+ \subset M$ is the conic isotropic subvariety
\begin{gather*}
\Lambda_+ = \Lambda \cup \coprod_{p\in \mathfrak c} \Lambda_{p+} \subset M.
\end{gather*}
\end{Definition}

\begin{Remark}
In the spirit of the adjunctions to come, one could note that we can rewrite the def\/inition of $\Lambda_{\tilde p} \subset C_p$
in the evidently equivalent forms
\begin{gather*}
\Lambda_{\tilde p} = q_p^{-1} q_p(\Lambda \cap C_p) = i_p q_p^{-1} q_p i_p^{-1}(\Lambda) \subset C_p.
\end{gather*}
\end{Remark}

\subsection{Open submanifolds}

We continue with $(M, \theta, \Lambda)$ a~marked Weinstein manifold.

Suppose that $p\in \mathfrak c$ is maximal, so that the coisotropic cell $i_p:C_p \to M$ is closed.
It will be useful to introduce some further structure on its normal geometry.

\begin{Definition}
Given a~marked Weinstein manifold $(M, \theta, \Lambda^\infty)$, and a~maximal critical point $p\in \mathfrak c$, a~{\em
defining function} $m_p:M \to [0, 1]$ for the corresponding closed coisotropic cell $C_p \subset M$, is a~subanalytic function
such that:
\begin{enumerate}\itemsep=0pt
\item[1)] $m_p^{-1}(0) = C_p$, \item[2)] $dm_p(Z) \leq 0$, \item[3)] $0, 1\in [0,1]$ are the only critical values of $m_p$, \item[4)] there is
an open subset $U\subset M$ containing $\mathfrak c \setminus \{p\} \subset M$ such that $U \subset m_p^{-1}(1)$, \item[5)] there is
a~compact subset $W \subset M$ such that $dm_p(Z) = 0$ over the complement $M \setminus W$.
\end{enumerate}
\end{Definition}

\begin{Lemma}
Defining functions always exist.
\end{Lemma}

\begin{proof}
This is easily obtained from the basic properties of subanalytic functions~\cite{BieMi, vdDM}.
\end{proof}

We will say that a~subset $\mathfrak s\subset \mathfrak c$ is open if $p\in \mathfrak s$ and $q\in \mathfrak c$ with $q\leq p$
implies $q\in \mathfrak s$.
Open subsets $\mathfrak s\subset \mathfrak c$ correspond to open unions of coisotropic cells
\begin{gather*}
M_\mathfrak s = \coprod_{p\in \mathfrak s} C_p.
\end{gather*}

\begin{Lemma}
\label{lemma open restriction}
For any open subset $\mathfrak s\subset \mathfrak c$, we obtain a~natural marked Weinstein manifold $(M_\mathfrak s,
\theta_\mathfrak s, \Lambda_\mathfrak s)$ by restriction
\begin{gather*}
\theta_\mathfrak s = \theta|_{M_\mathfrak s},
\qquad
\Lambda_\mathfrak s = \Lambda\cap M_\mathfrak s.
\end{gather*}
\end{Lemma}

\begin{proof}
We must produce a~Morse function $h_\mathfrak s:M_\mathfrak s\to {\mathbb R}$ such that $Z_\mathfrak s = Z|_{M_\mathfrak s} $ is
gradient-like with respect to $h_\mathfrak s$.
By induction, it suf\/f\/ices to assume that $\mathfrak c\setminus \mathfrak s$ is a~single maximal critical point $p$, so the
corresponding coisotropic cell $C_p \subset M$ is closed.

Let $h:M\to {\mathbb R}$ be a~Morse function such that $Z $ is gradient-like with respect to $h$.
Choose a~def\/ining function $m_p:M \to [0,1]$ for the closed coisotropic cell $C_p = M \setminus M_\mathfrak s$.
Consider the new function
\begin{gather*}
\xymatrix{ h_\mathfrak s = h + 1/m_p: M_\mathfrak s \ar[r] & {\mathbb R} }.
\end{gather*}
By construction, we have
\begin{gather*}
dh_\mathfrak s(Z_\mathfrak s) = dh(Z) - dm_p(Z)/m_p^2 \geq 0
\end{gather*}
with equality if and only if $dh = 0$ hence if and only if we are at a~point of $\mathfrak s$.
Furthermore, there is a~neighborhood of $\mathfrak s$ on which $dm_p = 0$, and so since $Z$ is gradient-like with respect to~$h$, we conclude $Z_\mathfrak s$ is gradient-like with respect to $h_\mathfrak s$.
\end{proof}

\begin{Remark}
In the case when $\Lambda^\infty = \varnothing$ so that $\Lambda = K$, we have $\Lambda_\mathfrak s = K \cap M_\mathfrak s$.
Note that the core $K_\mathfrak s \subset M_\mathfrak s$ is contained in $K \cap M_\mathfrak s$, with equality if and only if
$M_\mathfrak s = M$.
\end{Remark}

\section{Fukaya category}
\label{sec fukaya}

In this section, we survey the construction and basic properties of the stable inf\/initesimal Fukaya category of a~Weinstein
manifold~\cite{ N, NSpr, NZ, S01, Seidel}.

\subsection{Background structures}

Let $(M, \theta)$ be a~Weinstein manifold.
We brief\/ly review the standard additional data needed to consider Lagrangian branes and the Fukaya category.

We will work with compatible almost complex structures $J\in\operatorname{End}(TM)$ that are invariant under dilations near
inf\/inity.
The corresponding Riemannian metrics present $M$ near inf\/inity as a~metric cone over the projectivization $M^\infty$.
The space of all such compatible almost complex structures is nonempty and convex.

Given an almost complex structure $J\in\operatorname{End}(TM)$, we can speak about the complex canonical line bundle $\kappa_M =
(\wedge^{\dim M/2} T^{\rm hol} M)^{-1}$.
A bicanonical trivialization $\eta$ is an identif\/ication of the tensor-square $\kappa_M^{\otimes 2}$ with the trivial complex
line bundle.
The obstruction to a~bicanonical trivialization is twice the f\/irst Chern class $2c_1(M)\in H^2(M, {\mathbf Z})$, and all
bicanonical trivializations form a~torsor over the gauge group $\operatorname{Map}(M, S^1)$.
Forgetting the specif\/ic almost complex structure $J\in\operatorname{End}(TM)$, we will use the term bicanonical trivialization
to refer to a~section of bicanonical trivializations over the space of compatible almost complex structures.

Any symplectic manifold $M$ is canonically oriented, or in other words, after choosing a~Riemannian metric, the structure group
of $TM$ is the special orthogonal group.
A spin structure $\sigma$ is a~further lift of the structure group to the spin group.
The obstruction to a~spin structure is the second Stiefel--Whitney class $w_2(M)\in H^2(M, {\mathbf Z}/2{\mathbf Z})$, and all
spin structures form a~torsor over the group $H^1(M, {\mathbf Z}/2{\mathbf Z})$.

\begin{Definition}
A {\em Weinstein target} $(M, \theta, \eta, \sigma)$ is a~Weinstein manifold $(M, \theta)$ together with a~bicanonical
trivialization $\eta$ and spin structure $\sigma$.
We will often suppress mention of the latter structures when they are f\/ixed throughout.
\end{Definition}

\subsection{Lagrangian branes}

By an exact Lagrangian submanifold $L\subset M$, we mean a~closed but not necessarily compact submanifold of dimension $\dim
M/2$ such that the restriction $\theta|_{L}$ is an exact one-form, so in particular, $\omega|_{L} = 0$ where $\omega = d\theta$.

To ensure reasonable behavior near inf\/inity, we place two assumptions on our exact Lag\-rangian submanifolds $L \subset M$.
First, we insist that the compactif\/ication $\overline L \subset \overline M$ is a~subanalytic subset.
Along with other nice properties, this implies the following two facts:
\begin{enumerate}\itemsep=0pt
\item The boundary at inf\/inity $\partial^\infty L \subset M^\infty$ is an isotropic subvariety.
\item For $h:M\to {\mathbb R}$ a~Morse function such that $Z $ is gradient-like with respect to $h$, there is a~real number
$r>0$ such that the restricted function
\begin{gather*}
\xymatrix{ h: L \cap h^{-1}(r, \infty)\ar[r] & {\mathbb R} }
\end{gather*}
has no critical points.
\end{enumerate}
Second, we also assume the existence of a~perturbation $\psi$ that moves $L$ to a~nearby exact Lagrangian submanifold that is
tame (in the sense of~\cite{Sikorav}) with respect to a~conic metric.

\begin{Definition}
Fix a~f\/ield $k$.

A {\em brane structure} on an exact Lagrangian submanifold ${L} \subset M$ is a~three-tuple $ ({\mathcal E},\alpha,\flat) $
consisting of a~f\/lat f\/inite-dimensional $k$-vector bundle ${\mathcal E}\to L$, along with a~grading $\alpha$ (with respect to
the given compatible class of bicanonical trivializations) and a~pin structure $\flat$.

A {\em Lagrangian brane} in $M$ is a~four-tuple $ (L, {\mathcal E},\alpha,\flat) $ of an exact Lagrangian submanifold ${L}
\subset M$ equipped with a~brane structure $ ({\mathcal E},\alpha,\flat).
$ When there is no chance for confusion, we often write $L$ alone to signify the Lagrangian brane.

The {\em objects of the Fukaya category} $\operatorname{Ob} F(M)$ comprise all Lagrangian branes $L \subset M$.
\end{Definition}

Here is a~brief reminder on what a~grading and pin structure entail.
First, consider the bundle of Lagrangian planes ${\mathcal L}{\rm ag}_{M}\to M$, and the squared phase map
\begin{gather*}
\xymatrix{ \overline\alpha: {\mathcal L}{\rm ag}_{M} \ar[r] & {\mathbf C}^\times, & \overline\alpha(\mathcal L)= \eta(\wedge^{\dim M/2}
\mathcal L)^2 }.
\end{gather*}
Given a~Lagrangian submanifold $L \subset M$, we obtain the restricted map
\begin{gather*}
\xymatrix{ \overline\alpha:L\ar[r] & {\mathbf C}^\times, & \overline\alpha(x)=\overline\alpha(T_xL) }.
\end{gather*}
A grading of $L$ is a~lift
\begin{gather*}
\xymatrix{ \alpha:L \ar[r] & {\mathbf C}, & \overline \alpha = \exp\circ \alpha }.
\end{gather*}
The obstruction to a~grading is the Maslov class $[\overline\alpha]\in H^1(L,{\mathbf Z})$, and all gradings form a~torsor over
the group $H^0(L,{\mathbf Z})$.

Second, and recall that the (positive) pin group $\operatorname{Pin}^+(n)$ is the double cover of the orthogonal group $O(n)$ with center
${\mathbf Z}/2{\mathbf Z} \times {\mathbf Z}/2{\mathbf Z}$.
A pin structure $\flat$ on a~Riemannian manifold $L$ of dimension $n$ is a~lift of the structure group of $TL$ along the map
\begin{gather*}
\xymatrix{\operatorname{Pin}^+(n)\ar@{->>}[r] & O(n).
}
\end{gather*}
The obstruction to a~pin structure is the second Stiefel--Whitney class $w_2(L)\in H^2(L,{\mathbf Z}/2{\mathbf Z})$, and all
possible pin structures form a~torsor over the group $H^1(L,{\mathbf Z}/2{\mathbf Z})$.

\subsection{Intersections} Graded linear spans of intersection points provide the morphisms in the Fukaya category $F(M)$.
Given a~f\/inite collection of Lagrangians branes $L_0,\ldots, L_d\in \operatorname{Ob} F(M)$, we must perturb them so that their
intersections occur in some bounded domain.
To organize the perturbations, we recall the inductive notion of a~fringed set.

A fringed set $R_1\subset {\mathbb R}_+$ is any interval of the form $(0,r)$ for some $r>0$.
A fringed set $R_{d+1}\subset {\mathbb R}_+^{d+1}$ is a~subset satisfying the following:
\begin{enumerate}\itemsep=0pt
\item $R_{d+1}$ is open in $ {\mathbb R}^{d+1}_+$.
\item Under the projection $\pi: {\mathbb R}^{d+1}\to {\mathbb R}^d$ forgetting the last coordinate, the image $\pi(R_{d+1})$ is
a~fringed set.
\item If $(r_1,\ldots, r_d, r_{d+1})\in R_{d+1}$, then $(r_1,\ldots, r_d, r'_{d+1})\in R_{d+1}$ for $0<r'_{d+1}< r_{d+1}$.
\end{enumerate}

A Hamiltonian function $H:M\to {\mathbb R}$ is said to be {controlled} if near inf\/inity it is equal to a~conical coordinate.
Given a~f\/inite collection of Lagrangians branes $L_0,\ldots, L_d\in \operatorname{Ob} F(M)$, and controlled Hamiltonian
functions $H_0,\ldots, H_d$, we may choose a~fringed set $R\subset {\mathbb R}^{d+1}$ such that for $(\delta_d,\ldots,
\delta_0)\in R$, there is a~compact region $W\subset M$ such that for any $i\not = j$, we have
\begin{gather*}
\varphi_{H_i,\delta_i}(\overline L_i)\cap \varphi_{H_j,\delta_j}(\overline L_j)
\qquad
\mbox{lies in $W$.}
\end{gather*}
By a~further compactly supported Hamiltonian perturbation, we may also arrange so that the intersections are transverse.

We consider f\/inite collections of Lagrangian branes $L_0,\ldots, L_d\in \operatorname{Ob} F(M)$ to come equipped with such
perturbation data, with the brane structures $({\mathcal E}_i,\alpha_i, \flat_i)$ and taming perturbations $\psi_i$ transported
via the perturbations.

\begin{Definition}
Given a~f\/inite collection of Lagrangian branes $L_0,\ldots, L_d\in \operatorname{Ob} F(M)$, for branes $L_i$, $L_j$ with $i<j$, the
{\em graded vector space of morphisms} between them is the direct sum
\begin{gather*}
\hom_{F(M)}(L_i,L_j) = \bigoplus_{p\in \psi_i(\varphi_{H_i,\delta_i}(L_i))\cap \psi_j(\varphi_{H_j,\delta_j}(L_j)) } {\mathcal
H}{\rm om}(\mathcal E_i \vert_{p},\mathcal E_j\vert_{p})[-\deg(p)].
\end{gather*}
where the integer $\deg(p)$ denotes the Maslov grading at the intersection.
\end{Definition}

It is worth emphasizing that the salient aspect of the above perturbation procedure is the relative position of the perturbed
branes rather than their absolute position.
The following informal viewpoint can be a~useful mnemonic to keep the conventions straight.
In general, we always think of morphisms as ``propagating forward in time''.
Thus to calculate the morphisms $\hom_{F(M)}(L_0, L_1)$, we have required that $L_0$, $L_1$ are perturbed near inf\/inity so that~$L_1$ is further in the future than~$L_0$.
But what is important is not that they are both perturbed forward in time, only that $L_1$ is further along the timeline than~$L_0$.
So for example, we could perturb~$L_0$,~$L_1$ near inf\/inity in the opposite direction as long as $L_0$ is further in the past
than~$L_1$.

It is also worth noting the basic asymmetry of morphisms devolving from this def\/inition.
Because of the directionality of perturbations near inf\/inity, there is no general comparison of $\hom_{F(M)}(L_0, L_1)$ and
$\hom_{F(M)}(L_1, L_0)$ unlike when the branes are compact and a~Calabi--Yau relation holds.

\subsection{Compositions}

Signed counts of pseudoholomorphic polygons provide the dif\/ferential and higher composition maps of the $A_\infty$-structure of
the Fukaya category $F(M)$.
We use the following approach of Sikorav~\cite{Sikorav} (or equivalently, Audin--Lalonde--Polterovich~\cite{ALP}) to ensure that
the relevant moduli spaces are compact, and hence the corresponding counts are f\/inite.

First, a~Weinstein manifold $(M, \theta)$ equipped with a~compatible almost complex structure conic near inf\/inity is tame in the
sense of~\cite{Sikorav}.
To see this, it is easy to derive an upper bound on its curvature and a~positive lower bound on its injectivity radius.

Next, given a~f\/inite collection of branes $L_0,\ldots, L_d \in\operatorname{Ob} F(M)$, denote by $L$ the union of their
perturbations $\psi_i( \varphi_{H_i,\delta_i}(L_i))$ as described above.
By construction, the intersection of $L$ with the region $M \setminus W$ is a~tame submanifold (in the sense of~\cite{Sikorav}).
Namely, there exists $\rho_L>0$ such that for every $x\in L$, the set of points $y\in L$ of distance $d_{}(x,y) \leq \rho_L$ is
contractible, and there exists $C_L$ giving a~two-point distance condition $d_L(x,y) \leq C_L d_{}(x,y)$ whenever $x,y\in L$
with $d_{}(x,y)<\rho_L$.

Now, consider a~f\/ixed topological type of pseudoholomorphic map
\begin{gather*}
\xymatrix{ u:(D,\partial D) \ar[r] & (M, L).
}
\end{gather*}
Assume that all $u(D)$ intersect a~f\/ixed compact region, and there is an a~priori area bound ${\rm Area}(u(D))< A$.
Then as proven in~\cite{Sikorav}, one has compactness of the moduli space of such maps~$u$.
In fact, one has a~diameter bound (depending only on the given constants) constraining how far the image~$u(D)$ can stretch from
the compact set.

In the situation at hand, for a~given $A_\infty$-structure constant, we must consider pseudoholomorphic maps $u$ from polygons
with labeled boundary edges.
In particular, all such maps $u$ have image intersecting the compact set given by a~single intersection point.
The area of the image $u(D)$ can be expressed as the contour integral
\begin{gather*}
{\rm Area}(u(D)) = \int_{u(\partial D)} \theta.
\end{gather*}
Since each of the individual Lagrangian branes making up $L$ is exact, the contour integral only depends upon the integral of
$\theta$ along minimal paths between intersection points.
Thus such maps~$u$ satisfy an a~priori area bound.
We conclude that for each $A_\infty$-structure constant, the moduli space def\/ining the structure constant is compact, and its
points are represented by maps~$u$ with image bounded by a~f\/ixed distance from any of the intersection points.

\begin{Definition}
Given a~f\/inite collection of Lagrangian branes $L_0,\ldots, L_d\in \operatorname{Ob} F(M)$, the {\em composition map}
\begin{gather*}
m^d: {\hom}_{F(M)}(L_0,L_1)\otimes\cdots\otimes{\hom}_{F(M)}(L_{d-1},L_d)\rightarrow {\hom}_{F(M)}(L_0,L_d)[2-d]
\end{gather*}
is def\/ined as follows.
Consider elements $p_i\in {\hom}(L_i,L_{i+1})$, for $i=0,\ldots,d-1$, and $p_d\in \hom(L_0,L_d)$.
Then the coef\/f\/icient of $p_d$ in $m^d(p_0,\dots,p_{d-1})$ is def\/ined to be the signed sum over pseudoholomorphic maps from
a~disk with $d+1$ counterclockwise cyclically ordered marked points mapping to the~$p_i$ and corresponding boundary arcs mapping
to the perturbations of $L_{i+1}$.
Each map contributes according to the holonomy of its boundary, where adjacent perturbed components~$L_i$ and $L_{i+1}$ are
glued with~$p_i$.
\end{Definition}

\subsection{Coherence}

In the preceding sections, we have described the objects, morphisms, and compositions of the Fukaya category $F(M)$.
As explained in the fundamental sources~\cite{FOOO, Seidel}, there are a~large number of details to organize to be sure to
obtain an honest $A_\infty$-category.
In particular, calculations require branes be in general position, and hence must be invariant under suitable perturbations.
In the setting of noncompact branes, we have additional perturbations near inf\/inity to keep track of.
In particular, at f\/irst pass, the constructions given only provide what might be called a~compatible collection of directed
$A_n$-categories, for all $n$.
Here $A_n$ denotes Stashef\/f's operad of partial associative operations~\cite{Stasheff-I,Stasheff-II}, and we use the term directed as
in~\cite{Seidel}.
The former arises since we only def\/ine f\/initely many composition coef\/f\/icients at one time, and the latter since our
perturbations near inf\/inity are directed ``forward in time''.
To conf\/irm the coherence of the def\/initions, one can appeal to a~ref\/ined version of the well-known invariance of Floer
calculations under Hamiltonian isotopies.
We include a~brief discussion here (largely borrowed from~\cite{NSpr}) to explain the key ideas behind this approach.

Let $h:M\to \mathbb R$ be a~Morse function compatible with the Weinstein structure, in particular, providing a~conical
coordinate near inf\/inity.

\begin{Definition}
By a~{\em one-parameter family of closed $($but not necessarily compact$)$ submanifolds $($without boundary$)$} of~$M$, we mean a~closed
submanifold
\begin{gather*}
{\mathfrak L} \hookrightarrow \mathbb R\times M
\end{gather*}
satisfying the following:
\begin{enumerate}\itemsep=0pt
\item The restriction of the projection $p_\mathbb R: \mathbb R\times M\to \mathbb R$ to the submanifold ${\mathfrak L}$ is
nonsingular.

\item There is a~real number $r>0$, such that the restriction of the product $ p_\mathbb R\times h: \mathbb R\times M \to
\mathbb R \times \mathbb R $ to the subset $\{h>r\}\cap {\mathfrak L}$ is proper and nonsingular.

\item There is a~compact interval $[a,b]\hookrightarrow \mathbb R$ such that the restriction of the projection $p_M:\mathbb
R\times M\to M$ to the submanifold $p_\mathbb R^{-1}([\mathbb R\setminus [a,b])\cap {\mathfrak L}$ is locally constant.
\end{enumerate}
\end{Definition}

\begin{Remark}
Conditions (1) and (2) will be satisf\/ied if the restriction of the projection $\overline p_\mathbb R:\mathbb R\times \overline M
\to \mathbb R$ to the closure $\overline{\mathfrak L}\hookrightarrow \overline M$ is nonsingular as a~stratif\/ied map, but the
weaker condition stated is a~useful generalization.
It implies in particular that the f\/ibers ${\mathfrak L}_s = p_\mathbb R^{-1}(s) \cap{\mathfrak L} \hookrightarrow M$ are all
dif\/feomorphic, but imposes no requirement that their boundaries at inf\/inity should all be homeomorphic as well.
\end{Remark}

\begin{Definition}
By a~{\em one-parameter family of tame Lagrangian branes} in $M$, we mean a~one-parameter family of closed submanifolds
${\mathfrak L}\hookrightarrow \mathbb R\times M$ such that the f\/ibers ${\mathfrak L}_s = p_\mathbb R^{-1}(s) \cap{\mathfrak L}
\hookrightarrow M$ also satisfy:
\begin{enumerate}\itemsep=0pt

\item The f\/ibers ${\mathfrak L}_s$ are exact tame Lagrangians with respect to the symplectic structure and any almost complex
structure conical near inf\/inity.

\item The f\/ibers ${\mathfrak L}_s$ are equipped with a~locally constant brane structure $(\mathcal E_s,\alpha_s,\flat_s)$ with
respect to the given background classes.

\end{enumerate}
\end{Definition}

\begin{Remark}
Note that if we assume that ${\mathfrak L}_0$ is an exact Lagrangian, then ${\mathfrak L}_s$ being an exact Lagrangian is equivalent
to the family ${\mathfrak L}$ being given by the f\/low $\varphi_{H_s}$ of the vector f\/ield of a~time-dependent Hamiltonian
$H_s:M\to \mathbb R$.
Note as well that a~brane structure consists of topological data, so can be transported unambiguously along the f\/ibers of such
a~family.
\end{Remark}

\begin{Remark}
It is important for various applications of the theory of this section that we consider one-parameter families of tame
Lagrangian branes ${\mathfrak L}\hookrightarrow \mathbb R\times M$ such that the f\/iber branes ${\mathfrak L}_s \subset M$ are not
constant near inf\/inity.
Thus in general the f\/iber branes will not be isomorphic objects of $F(M)$, though we will see that their Floer-pairing with
respect to certain test branes will be invariant.
\end{Remark}

The rest of this section will be devoted to the following statement of Floer invariance.
It is the basic instance (going beyond the foundational results of~\cite{FOOO, Seidel}) of the general pattern that conf\/irms
$F(M)$ is a~well-def\/ined $A_\infty$-category.

\begin{Proposition}
\label{prop inv}
Suppose ${\mathfrak L}_s$ is a~one-parameter family of tame Lagrangian branes in $M$.
Suppose $L'$ is a~fixed test brane which is disjoint from ${\mathfrak L}_s$ near infinity for all $s$.
Suppose ${\mathfrak L}_s$ is transverse to $L'$ except for finitely many points.

Then for any $a$, $b$ with ${\mathfrak L}_a$ and ${\mathfrak L}_b$ transverse to $L'$, the Floer chain complexes ${\rm CF}({\mathfrak L}_{a},
L')$ and ${\rm CF}({\mathfrak L}_{b}, L')$ are quasi-isomorphic.
\end{Proposition}

Before proving the proposition in general, it is convenient to f\/irst prove the following special case.

\begin{Lemma}
Suppose ${\mathfrak L}_s$ is one-parameter family of tame Lagrangian branes in~$M$.
Suppose~$L'$ is a~f\/ixed test object which is disjoint from ${\mathfrak L}_s$ near inf\/inity for all~$s$.

Fix $s_0$ and assume ${\mathfrak L}_{s_0}$ is transverse to $L'$.
Then there is an $\epsilon >0$ so that for all $s_1 \in (s_0 - \epsilon, s_0+\epsilon)$,
the Floer chain complexes ${\rm CF}({\mathfrak L}_{s_0}, L')$ and ${\rm CF}({\mathfrak L}_{s_1}, L')$ are quasi-isomorphic.
\end{Lemma}

\begin{proof}
By our assumptions on the tame behavior (in the sense of~\cite{Sikorav}) of ${\mathfrak L}_{s_0}$ and $L'$ near inf\/inity, the
moduli spaces giving the dif\/ferential of ${\rm CF}({\mathfrak L}_{s_0}, L')$ are compact.
This follows from the {a priori} $C^0$-bound: there is some $r_0\gg 0$, such that no disk in the moduli space leaves the region
$h < r_0$, where $h:M\to \mathbb R$ is a~compatible Morse function.

Choose some $r_1 > r_0$.
Then for very small $\epsilon> 0$ and any $s_1 \in (s_0 - \epsilon, s_0+\epsilon)$,
we may decompose the motion ${\mathfrak L}_{s_0} \rightsquigarrow {\mathfrak L}_{s_1}$
into two parts: f\/irst, a~motion ${\mathfrak L}_{s_0} \rightsquigarrow L$ {\em
supported in the region $h > r_1$}; and then second, a~{\em compactly supported} motion $L \rightsquigarrow {\mathfrak L}_{s_1}$.
We must show that each of the above two motions leads to a~quasi-isomorphism.

First, for the motion ${\mathfrak L}_{s_0} \rightsquigarrow L$, since we have not changed ${\mathfrak L}_{s_0}$ or $L'$ in the
region $h < r_0$, the same {a priori} $C^0$-bounds of~\cite{Sikorav} hold (they only depend on the Lagrangians in the region $h
< r_0$), and the pseudoholomorphic strips for the pair $({\mathfrak L}_{s_0}, L')$ and for the pair $(L, L')$ are in fact {exactly
the same} (we could perversely attach ``wild'' non-intersecting ends to either and it would not make a~dif\/ference.) Thus we can
take the continuation map to be the identity.

Second, the motion $L \rightsquigarrow {\mathfrak L}_{s_1}$ is compactly supported, so standard PDE techniques provide
a~continuation map.
\end{proof}

\begin{Remark}
In the proof of the lemma, one should probably not use the term ``continuation map'' for the quasi-isomorphism associated to the
f\/irst motion ${\mathfrak L}_{s_0} \rightsquigarrow L$.
Rather, it is an example of the more general setup of parameterized moduli spaces.
In the above setting, one can obtain a~uniform $C^0$-bound over the family, so the parameterized moduli space is compact, and
hence one can apply standard cobordism arguments to prove the matrix coef\/f\/icients at the initial and f\/inal time are the same.
\end{Remark}

\begin{proof}[Proof of Proposition~\ref{prop inv}] {\sloppy By the previous lemma, it suf\/f\/ices to show that for any~$s_0$ with~${\mathfrak L}_{s_0}$ not
(necessarily) transverse to~$L'$, there is a~small $\epsilon>0$ such that the Floer chain complexes ${\rm CF}({\mathfrak L}_{s_0
-\epsilon}, L')$ and ${\rm CF}({\mathfrak L}_{s_0+\epsilon}, L')$ are quasi-isomorphic.

}

To see this, f\/ix a~compatible Morse function $h:M\to\mathbb R$, and a~bump function $b:\mathbb R\to\mathbb R$ such that the
composition $b\circ h:M\to \mathbb R$ is identically $0$ near inf\/inity and $1$ on a~compact set containing all of the critical
points of $h$ and (possibly non-transverse) intersection points ${\mathfrak L}_{s_0} \cap L'$.

Let $H_s:M\to \mathbb R$ be a~(time-dependent) Hamiltonian giving the motion ${\mathfrak L}_s$.
The product Hamiltonian $\tilde H_s = (b\circ h) \cdot H_s:M\to \mathbb R$ gives a~family $\tilde {\mathfrak L}_s$ through the
base object ${\mathfrak L}_{s_0}$ satisfying: (1) $\tilde {\mathfrak L}_s$ is transverse to $L'$ whenever $|s-s_0| $ is small and
nonzero, and (2) $\tilde {\mathfrak L}_s$ is equal to ${\mathfrak L}_{s_0}$ near inf\/inity.
Therefore since the motion of $\tilde {\mathfrak L}_s$ is compactly supported, standard PDE techniques provide a~continuation map
giving a~quasi-isomorphism between ${\rm CF}(\tilde {\mathfrak L}_{s_0 -\epsilon}, L')$ and ${\rm CF}(\tilde {\mathfrak L}_{s_0 + \epsilon},
L')$, for small enough $\epsilon >0$.

Finally, returning to the bump function $b$, one can construct motions ${\mathfrak L}_{s_0 - \epsilon} \rightsquigarrow \tilde
{\mathfrak L}_{s_0 - \epsilon}$ and $\tilde {\mathfrak L}_{s_0 +\epsilon} \rightsquigarrow {\mathfrak L}_{s_0 + \epsilon}$ which are
supported near inf\/inity and thus in particular always transverse to $L'$.
Thus we may apply the previous lemma to obtain quasi-isomorphisms between ${\rm CF}({\mathfrak L}_{s_0 -\epsilon}, L')$ and ${\rm CF}( \tilde
{\mathfrak L}_{s_0 - \epsilon}, L')$, and similarly, between ${\rm CF}(\tilde {\mathfrak L}_{s_0 +\epsilon}, L')$ and ${\rm CF}( {\mathfrak L}_{s_0 + \epsilon}, L')$.
Putting together the above, we obtain a~quasi-isomorphism between ${\rm CF}({\mathfrak L}_{s_0 -\epsilon}, L')$ and ${\rm CF}( {\mathfrak L}_{s_0 + \epsilon}, L')$.
\end{proof}

\begin{Remark}
The above proposition (which is a~condensed form of arguments of~\cite{N, NZ} and appears explicitly in~\cite{NSpr} for
cotangent bundles) is closely related to Question~1.3 of Oh's paper~\cite{Oh} which asks whether a~homology-level continuation
map constructed by a~careful limiting argument with PDE techniques is induced by a~chain-level morphism.
While we have not investigated this, it is not hard to believe that the quasi-isomorphism of the above proposition provides the
desired lift.
\end{Remark}

\subsection{Stabilization}

It is convenient to work interchangeably with small idempotent-complete pre-triangulated $k$-linear
$A_\infty$-categories~\cite{FOOO, Seidel} when thinking about Fukaya categories and small stable idempo\-tent-complete $k$-linear
$\infty$-categories~\cite{HA, topos} when thinking about abstract constructions.
They have equivalent homotopy theories, and we lose nothing by going back and forth.
In what follows, all of the specif\/ic assertions we will use can be found in~\cite{BFN}.

\begin{Definition}
Let ${\mathcal C}$, ${\mathcal C}'$ be stable $\infty$-categories.
A functor $F:{\mathcal C}\to {\mathcal C}'$ is said to be
\begin{enumerate}\itemsep=0pt
\item[(1)] {\em continuous} if it preserves coproducts,

\item[(2)] {\em proper} if it preserves compact objects,

\item[(3)] {\em exact} if it preserves zero objects and exact triangles (equivalently, f\/inite colimits).
\end{enumerate}
\end{Definition}

It is convenient to work alternatively within two related $k$-linear contexts.

\begin{Definition}
We denote by $\operatorname{St}_k$ the symmetric monoidal $\infty$-category of stable presentable $k$-linear $\infty$-categories with morphisms
continuous functors.
The monoidal unit is the $\infty$-category $\operatorname{Mod} k$ of $k$-chain complexes.

We denote by $\operatorname{st}_k$ the symmetric monoidal $\infty$-category of small stable idempotent-complete $k$-linear $\infty$-categories
with morphisms exact functors.
The monoidal unit is the $\infty$-category $\operatorname{Perf} k$ of perfect $k$-chain complexes.
\end{Definition}

\begin{Definition}

The {\em big stabilization} $\operatorname{Mod} {\mathcal C} \in \operatorname{St}_k$ of a~small $k$-linear $\infty$-category ${\mathcal C}$
is the stable presentable $k$-linear dif\/ferential graded category of $A_\infty$-right modules
\begin{gather*}
\xymatrix{ {\mathcal C}^{\operatorname{op}} \ar[r] & \operatorname{Mod} k.
}
\end{gather*}

The {\em small stabilization} $\operatorname{Perf} {\mathcal C}\in \operatorname{st}_k $ is the small $k$-linear full $\infty$-subcategory of
$\operatorname{Mod} {\mathcal C}$ consisting of compact objects (summands of f\/inite colimits of representable objects).
\end{Definition}

\begin{Lemma}
For a~small $k$-linear $\infty$-category ${\mathcal C}$, the {Yoneda embedding} is fully faithful
\begin{gather*}
\xymatrix{ {\mathcal Y}:{\mathcal C} \ar@{^(->}[r] & \operatorname{Perf} {\mathcal C}, & {\mathcal Y}_L(P) = \hom_{{\mathcal
C}}(P, L) }.
\end{gather*}

If ${\mathcal C}\in \operatorname{st}_k$, then the Yoneda embedding is an equivalence.
\end{Lemma}

\begin{Corollary}
For a~small $k$-linear $\infty$-category ${\mathcal C}$, forming its stabilization canonically commutes with forming its
opposite category
\begin{gather*}
\operatorname{Perf}({\mathcal C}^{\operatorname{op}}) \simeq \operatorname{Perf}({\mathcal C})^{\operatorname{op}}.
\end{gather*}

\end{Corollary}

\begin{Lemma}
Forming big stabilizations is a~faithful symmetric monoidal functor
\begin{gather*}
\xymatrix{ \operatorname{Mod} :\operatorname{st}_k \ar[r]& \operatorname{St}_k }.
\end{gather*}

The monoidal dual of ${\mathcal C} \in \operatorname{st}_k$ is the opposite category ${\mathcal C}^{\operatorname{op}}$.

The monoidal dual of $\operatorname{Mod} {\mathcal C} \in \operatorname{St}_k$ is the restricted opposite category $\operatorname{Mod}
({\mathcal C}^{\operatorname{op}})$.
\end{Lemma}

\begin{Remark}
We can recover ${\mathcal C}\in \operatorname{st}_k$ from $\operatorname{Mod} {\mathcal C}\in \operatorname{St}_k$ by passing to compact objects ${\mathcal
C} = (\operatorname{Mod} {\mathcal C})^c$.
The image of the morphism $\operatorname{Mod} :\hom_{\operatorname{st}_k}({\mathcal C}, {\mathcal C}')\to \hom_{\operatorname{St}_k}(\operatorname{Mod}
{\mathcal C}, \operatorname{Mod} {\mathcal C}')$ comprises proper functors.
\end{Remark}

\begin{Corollary}
\label{cor hom vs tensor}
For ${\mathcal C}, {\mathcal C}' \in \operatorname{st}_k$, there are canonical equivalences
\begin{gather*}
\xymatrix{ {\mathcal C}^{\operatorname{op}} \otimes {\mathcal C}' \ar[r]^-\sim & \hom_{\operatorname{st}_k}({\mathcal C}, {\mathcal C}') },
\\
\xymatrix{ \operatorname{Mod} ({\mathcal C}^{\operatorname{op}} \otimes {\mathcal C}') \simeq \operatorname{Mod} ({\mathcal C}^{\operatorname{op}})
\otimes \operatorname{Mod} ({\mathcal C}') \ar[r]^-\sim & \hom_{\operatorname{St}_k}(\operatorname{Mod} {\mathcal C}, \operatorname{Mod} {\mathcal C}') }.
\end{gather*}
\end{Corollary}

\begin{Definition}
The {\em perfect Fukaya category} $\operatorname{Perf} F(M)$ is the small stabilization of $F(M)$.

The {\em stable Fukaya category} $\operatorname{Mod} F(M)$ is the big stabilization of $F(M)$.
\end{Definition}

\subsection{Singular support}
\label{sing supp}

While calculations among Lagrangian branes ref\/lect quantum topology, we nevertheless have access to their underlying Lagrangian
submanifolds.
We take advantage of this in the following def\/inition.

\begin{Definition}
Fix an object $L\in \operatorname{Perf} F(M)$.
\begin{enumerate}\itemsep=0pt
\item[(1)] The {\em null locus} $n(L) \subset M$ is the conic open subset of points $x\in M$ for which there exists a~conic open set $U
\subset M$ containing $x$ such that we have the vanishing
\begin{gather*}
\hom_{F(M)}(L,P)\simeq 0,
\qquad
\text{for any}
\quad
P\in F(M)
\quad
\text{with}
\quad
P \subset U, \quad \partial^\infty P\subset U^\infty.
\end{gather*}

\item[(2)] The {\em singular support} $\operatorname{ss}(L )\subset M$ is the conic closed complement
\begin{gather*}
\operatorname{ss}(L) = M \setminus n(L).
\end{gather*}
\end{enumerate}
\end{Definition}

\begin{Remark}
For cotangent bundles, under the equivalence of branes and constructible sheaves recalled in the next section, the above notion
of singular support for branes coincides with the traditional notion of singular support of sheaves.
\end{Remark}

\begin{Definition}
Let $(M, \theta, \Lambda)$ be a~marked Weinstein manifold.

We def\/ine the full subcategory $\operatorname{Perf}_\Lambda F(M) \subset \operatorname{Perf} F(M)$ to comprise objects $L\in
\operatorname{Perf} F(M)$ with singular support satisfying $\operatorname{ss}(L) \subset \Lambda$.
\end{Definition}

\begin{Lemma}
$\operatorname{Perf}_\Lambda F(M)$ is a~small stable idempotent-complete $k$-linear $\infty$-category.

If $\Lambda_1 \subset \Lambda_2$, then $\operatorname{Perf}_{\Lambda _1}F(M)\subset \operatorname{Perf}_{\Lambda_2} F(M)$, and
$\operatorname{Perf} F(M) = \cup_{\Lambda} \operatorname{Perf}_\Lambda F(M)$.
\end{Lemma}

\begin{proof}
The singular support condition is clearly preserved by extensions and summands.
\end{proof}

In the remainder of this section, we explain (without proof) how to calculate the projectivization of the singular support.
By induction using recollement, this provides a~complete picture of the singular support.
We will not need this material for any further developments, but include it to help orient the reader.

Fix a~f\/inite collection of Lagrangian branes $L_i\in F(M)$, for $i \in I$, and let $V = \cup_{i \in I} L_i \subset M$ denote the
Lagrangian subvariety given by their union.

\begin{Lemma}
The boundary at infinity $\partial^\infty V \subset M^\infty$ is a~closed Legendrian subvariety.
\end{Lemma}

\begin{proof}
By dilation, we can contract $V\subset M$ to a~conical Lagrangian subvariety $V^{c} \subset M$.
Then we need only observe that $\partial^\infty V = \partial^\infty V^{c}$.
\end{proof}

Let $\partial_{sm}^\infty V \subset \partial^\infty V$ denote the smooth locus.
Given a~point $x \in \partial_{sm}^\infty V$, we can f\/ind a~small Legendrian sphere $S(x) \subset M^\infty$ centered at $x$, and
simply linked around $\partial_{sm}^\infty V$.
Then we can f\/ind a~Lagrangian brane $B(x)\subset M$ dif\/feomorphic to a~ball, and with boundary at inf\/inity $\partial^\infty B(x)
= S(x)$.
The particular grading on $B(x)$ will play no role.

\begin{Proposition}
For an object $L\in \operatorname{Perf} F(M)$ in the perfect envelope of the finite collection $L_i\in F(M)$, for $i \in I$, the
projectivization of its singular support $\operatorname{ss}( L)^\infty \subset M^\infty$ is the closure of the subset
\begin{gather*}
\big\{x \in \partial_{sm}^\infty V \,|\, \hom_{F(M)} (L, B(x)) \not \simeq 0\big\}.
\end{gather*}
\end{Proposition}

\begin{Corollary}
The projectivization of the singular support $\operatorname{ss}(L)^\infty \subset M^\infty$ is a~Legendrian subvariety.
\end{Corollary}

\subsection{Cotangent bundles} We brief\/ly remind the reader of the equivalence between branes in a~cotangent bundle and
constructible sheaves on the base manifold.

Let $X$ be a~compact manifold with cotangent bundle $\pi:T^*X \to X$ and projectivization $\pi^\infty:T^\infty X \to X$.
For simplicity, let us assume that $X$ is equipped with an orientation and spin structure.
Then as explained for instance in Example~\ref{cotangent is weinstein}, we can view $T^*X$ as a~Weinstein manifold with
a~canonically trivial canonical bundle and canonical spin structure.

\begin{Definition}
Let $\operatorname{Sh}(X)$ denote the dif\/ferential graded category of complexes of sheaves on~$X$ with constructible cohomology.
\end{Definition}

\begin{Definition}
Given an object ${\mathcal F}\in \operatorname{Sh}(X)$, we write $\operatorname{ss} {\mathcal F}\subset T^*X$ for its singular support, and
$\operatorname{ss}^\infty {\mathcal F} \subset T^\infty X$ for the projectivization of its singular support.

Given a~Whitney stratif\/ication ${\mathcal S} =\{X_\alpha\}_{\alpha\in A}$, we def\/ine its conormal bundle and projectivized
conormal bundle to be the unions
\begin{gather*}
T^*_{\mathcal S} X = \coprod_{\alpha \in A} T^*_{X_\alpha} X \subset T^* X,
\qquad
T^\infty_{\mathcal S} X=\coprod_{\alpha \in A} T^\infty_{X_\alpha} X \subset T^\infty X.
\end{gather*}
\end{Definition}

\begin{Definition}
Suppose $\Lambda\subset T^*X$ is a~conical Lagrangian subvariety.
Let $\operatorname{Sh}_\Lambda(X)\subset \operatorname{Sh}(X)$ denote the full subcategory of complexes of sheaves with $\operatorname{ss}
{\mathcal F} \subset \Lambda$.

Suppose ${\mathcal S} =\{X_\alpha\}_{\alpha\in A}$ is a~Whitney stratif\/ication of $X$.
Let $\operatorname{Sh}_{\mathcal S}(X)\subset \operatorname{Sh}(X)$ denote the full subcategory of complexes of sheaves with
${\mathcal S}$-constructible cohomology.
\end{Definition}

\begin{Lemma}
For a~Whitney stratification ${\mathcal S}$, we have $\operatorname{Sh}_{{\mathcal S}} (X) = \operatorname{Sh}_{T^*_{\mathcal S}
X}(X)$.
\end{Lemma}

Now let $i:Y\to X$ be a~locally closed submanifold with frontier $\partial Y = \overline Y \setminus Y$.
On the one hand, we have the standard and costandard extensions $i_*k_{Y}, i_!k_{Y }\in \operatorname{Sh}(X)$.

On the other hand, we have corresponding Lagrangian branes constructed as follows.
Choose a~non-negative function $m:X\to{\mathbb R}_{\geq 0}$ with zero-set precisely $\partial Y \subset X$.

\begin{Definition}
We def\/ine the {\em standard} and {\em costandard Lagrangians} $L_{Y*}, L_{Y!}\in F(T^*X)$ to be the f\/iberwise translations
\begin{gather*}
L_{Y*}= \Gamma_{d\log m} + T^*_Y X,
\qquad
L_{Y!}= -\Gamma_{d\log m} + T^*_Y X
\end{gather*}
equipped with the orientation bundle $\operatorname{or}_{Y}$, and canonical gradings and spin structures.
\end{Definition}

We have the following from~\cite{ N, NSpr, NZ}.

\begin{Theorem}
\label{micro}
There is a~canonical equivalence
\begin{gather*}
\xymatrix{ \mu_X:\operatorname{Sh}(X) \ar[r]^-\sim & \operatorname{Perf} F(T^* X) }
\end{gather*}
such that $\mu_X(i_*k_Y) \simeq L_{Y*}$, $\mu_X(i_!k_Y) \simeq L_{Y!}$.

Furthermore, we have $\operatorname{ss}^\infty {\mathcal F} = \operatorname{ss}^\infty \mu_X({\mathcal F})$, and hence for a~conical Lagrangian subvariety
$\Lambda\subset T^*X$ containing the zero section, $\mu_X$ restricts to an equivalencez
\begin{gather*}
\xymatrix{ \mu_X:\operatorname{Sh}_\Lambda(X) \ar[r]^-\sim & \operatorname{Perf}_\Lambda F(T^* X) }.
\end{gather*}
\end{Theorem}

There are various extensions of the above result to noncompact manifolds $X$, but we will only call upon the following expected
generalization in support of the anticipated recollement pattern sketched in the next section.

Suppose $X$ is a~manifold whose noncompactness is concentrated at single conical end.
In other words, we have a~manifold $Y$ with boundary $\partial Y$ such that $\overline Y = Y \coprod \partial Y$ is compact, and
an identif\/ication
\begin{gather*}
X \simeq \overline Y \coprod_{\partial Y} (\partial Y \times [0, \infty)).
\end{gather*}
We can equivalently assume $X$ is equipped with a~Morse function $f_X:X\to {\mathbb R}$ and Riemannian metric $g$ such that the
f\/low of the gradient $\nabla_g f_X$ exhibits the conical end as a~metric product.

Consider the f\/iberwise linear function $F_X:T^*X \to {\mathbb R}$ def\/ined by $F_X(x, \xi) = \xi(\nabla_g f_X|_x)$, and the Morse
function $h = g + \pi^*_Xf_X: T^*X\to {\mathbb R}$ def\/ined by $h(x, \xi) = |\xi|_g^2 + f_X(x)$.
For $\epsilon>0$ suf\/f\/iciently small, the pair $(T^*X, \theta_X + \epsilon dF_X)$ should form a~Weinstein manifold exhibited by
the Morse function $h$.

Now suppose in addition that ${\mathcal S} =\{X_\alpha\}_{\alpha\in A}$ is a~Whitney stratif\/ication with strata $X_\alpha
\subset X$ that are conical near the end.
In other words, we assume that the strata $X_\alpha$ are transverse to~$\partial Y$ inside of~$X$, and the above identif\/ication
restricts to an identif\/ication
\begin{gather*}
X_\alpha = (X_\alpha \cap \overline Y) \coprod_{(X_\alpha \cap \partial Y)} ((X_\alpha \cap \partial Y)\times [0,\infty)).
\end{gather*}
We can equivalently assume the strata are invariant under the f\/low of the gradient of our Morse function along the conical end.

Now let $i:Y\to X$ be a~union of strata.
On the one hand, we have the standard and costandard extensions $i_*k_{Y}, i_!k_{Y }\in \operatorname{Sh}_{\mathcal S}(X)$.

On the other hand, we have the corresponding standard and costandard branes $L_{Y *}, L_{Y !}\in F(T^*X)$ constructed with
a~non-negative function $m:X\to{\mathbb R}_{\geq 0}$ invariant under the f\/low of the gradient of our Morse function along the
conical end.

With only slight modif\/ications, the arguments used to establish Theorem~\ref{micro} should provide the following proposition.
To realize it, one could embed $X$ as an open submanifold of a~compact manifold $Z$ and interpret all of the calculations to
take place there.
In particular, one could extend all sheaves by zero of\/f of $X$ and work with the corresponding branes in $T^*Z$ as prescribed by
Theorem~\ref{micro}.
Thanks to our perturbation framework and diameter bounds on disks, the relevant geometry should be equivalent to that within
$T^*X$.

Now we expect there is a~canonical functor
\begin{gather*}
\xymatrix{ \mu_X:\operatorname{Sh}_{\mathcal S}(X) \ar[r]& \operatorname{Perf} F(T^* X) }
\end{gather*}
such that $\mu_X(i_{*}k_{Y}) \simeq L_{Y*}$, $\mu_X(i_{!}k_{Y}) \simeq L_{Y!}$.

\begin{Remark}
We will return to the above in the special situation when $X$ itself is a~Weinstein manifold, with $f_X$ a~Morse function
compatible with its Liouville vector f\/ield, and ${\mathcal S}$ the stratif\/ication by coisotropic cells.
\end{Remark}

\section{Adjunctions for branes}

Due to the technical demands of the material, this section will be more informal than the others.
We outline the expected recollement pattern of diagrams~\eqref{eq functors} and~\eqref{eq triangles}, only sketching the
constructions and results that we anticipate should hold.

\subsection{Bimodules via correspondences}

\begin{Definition}
Let $(M, \theta)$ be an exact symplectic manifold.
The {\em opposite} exact symplectic manifold $(M^{\operatorname{op}}, -\theta)$ is the same underlying manifold $M$ equipped with the negative
Liouville form, and hence negative symplectic form $-\omega = -d\theta$.
\end{Definition}

When $(M, \theta)$ is equipped with background structures, we transport them by the identity to obtain background structures on
$(M^{\operatorname{op}}, -\theta)$.
In particular, given a~compatible almost complex structure $J\in\operatorname{End}(TM)$, we take $-J\in
\operatorname{End}(TM^{\operatorname{op}})$.
Thus we have an identity of bicanonical line bundles $\kappa_{M^{\operatorname{op}}} \simeq \kappa_{M}^{-1}$ and so a~bicanonical
trivialization for $(M, \theta)$ induces one for $(M^{\operatorname{op}}, -\theta)$.

\begin{Lemma}
If $(M, \theta)$ is a~Weinstein manifold, its opposite $(M^{\operatorname{op}}, -\theta)$ is a~Weinstein manifold with the same Liouville
vector field.
\end{Lemma}

\begin{proof}
If $\theta = i_{Z}\omega$, then $-\theta = i_{Z}(-\omega)$.
\end{proof}

\begin{Proposition}
There is a~canonical identification
\begin{gather*}
\xymatrix{ F(M) \ar[r]^-\sim & F(M^{\operatorname{op}})^{\operatorname{op}} }
\end{gather*}
given on Lagrangian branes by the duality
\begin{gather*}
\xymatrix{ L = (L, {\mathcal E}, \alpha, \flat) \ar@{|->}[r] & L' = (L, {\mathcal E}^\vee, -\alpha, \flat) }.
\end{gather*}
\end{Proposition}

\begin{proof}
Our perturbation framework is compatible with the assertion.
\end{proof}

\begin{Definition}
Suppose $M$, $N$ are Weinstein targets.

A {\em Lagrangian correspondence} is an object $\mathsf K\in F(M^{\operatorname{op}} \times N)$.
The {\em dual correspondence} $\mathsf K'\in F(N^{\operatorname{op}}\times M )$
is the matched object under the equivalence $F(M^{\operatorname{op}} \times N)
\simeq F(N^{\operatorname{op}} \times M)^{\operatorname{op}}$.
\end{Definition}

Suppose $M, N$ are Weinstein targets, and $\mathsf K\in F(M^{\operatorname{op}} \times N)$ a~Lagrangian correspondence.
Let $L\in F(M)$, $P\in F(N)$ be test branes.

One should be able to obtain a~functor
\begin{gather*}
\xymatrix{ f_\mathsf K: F(M)\ar[r] & \operatorname{Mod} F(N) }
\end{gather*}
or equivalently a~morphism
\begin{gather*}
\xymatrix{ f_\mathsf K \in \hom_{\operatorname{St}_k}(\operatorname{Mod} F(M), \operatorname{Mod} F(N)) \simeq \operatorname{Mod} (F(M^{\operatorname{op}})
\otimes F(N)) }
\end{gather*}
given by the functorial construction
\begin{gather*}
\xymatrix{ f_\mathsf K(L)(P) = \hom_{F(M^{\operatorname{op}} \times N)}(L' \times P, \mathsf K).
}
\end{gather*}

There are (at least) two approaches one might take to conf\/irm that the above def\/inition is sensible.
The two approaches we will mention here should lead to homotopically equivalent theories.

On the one hand, one could adopt the geometric formalism of Lagrangian correspondences as developed by Wehrheim and
Woodward~\cite{ww6,ww3, ww1,ww4,ww2, ww5} and Ma'u~\cite{Mau} and count pseudoholomorphic quilts to provide the structure
constants of an $A_\infty$-functor.

On the other hand, one should be able to adopt the algebraic formalism of bimodules.
To implement this, unwinding the Yoneda embedding, one needs to know that the product on branes induces a~bilinear functor
\begin{gather*}
\xymatrix{ F(M^{\operatorname{op}}) \times F(N) \ar[r] & F(M^{\operatorname{op}} \times N) }.
\end{gather*}
With this in hand, one could then def\/ine the module $f_\mathsf K(L)$ as the composition
\begin{gather*}
\xymatrix{ f_\mathsf K(L):F(N) \ar[r] & \operatorname{Mod} k, & f_\mathsf K(L)(P) = \hom_{F(M^{\operatorname{op}} \times N)}(L' \times P,\mathsf K).}
\end{gather*}
This should have the simultaneous functoriality in $P$ to make $f_K(L)$ a~module, and the functoriality in $L$ to make $f_K$
a~functor.
Note that establishing the existence of the product functor is much less than proving that the induced linear functor
\begin{gather*}
\xymatrix{ F(M^{\operatorname{op}}) \otimes F(N) \ar[r] & F(M^{\operatorname{op}} \times N) }
\end{gather*}
is an equivalence.
Coupled with Corollary~\ref{cor hom vs tensor} in mind, this stronger assertion would show all functors are given by kernels,
but this is not necessary for the aims outlined above.
As far as we know, the existence of the product functor is not in the literature, though related homotopical algebra is
available~\cite{loday, MS, SU1, SU2} and some readers may prefer this point of view to the intricacies of Lagrangian
correspondences and pseudoholomorphic quilts.

{\sloppy With the appropriate brane structure, the diagonal $\Delta_M \subset M^{\operatorname{op}} \times M$ Lagrangian corresponden\-ce should give an
endofunctor $f_{\Delta_M}$ of $F(M)$ canonically equivalent to the iden\-ti\-ty~$\operatorname{id}_{F(M)}$.
A~si\-milar assertion in parallel settings can be found in many places inclu\-ding~\cite{ww6,ww3, ww1,ww4,ww2, ww5}.

}

\begin{Remark}
The above integral transform formalism is tuned to constructing right adjoints.
There is an alternative dual formalism suited to constructing left adjoints.

Given a~Lagrangian correspondence $\mathsf K\in F(M^{\operatorname{op}} \times N)$, and test branes $L\in F(M)$, $P\in F(N)$, one should be able
also to obtain a~functor
\begin{gather*}
\xymatrix{ g_\mathsf K: F(M)^{\operatorname{op}} \ar[r] & \operatorname{Mod} (F(N)^{\operatorname{op}}) }
\end{gather*}
given by the functorial construction
\begin{gather*}
\xymatrix{ g_\mathsf K(L)(P) = \hom_{F(M^{\operatorname{op}} \times N)}(\mathsf K, L' \times P).
}
\end{gather*}
Note this should be contravariant in the Lagrangian correspondence $\mathsf K$.

If the above produces a~proper functor, then in turn it can be interpreted as a~functor $g_\mathsf K:\operatorname{Perf} F(M)\to
\operatorname{Perf} F(N)$.

With the appropriate brane structure, the diagonal $\Delta_M \subset M^{\operatorname{op}} \times M$ Lagrangian correspondence should give an
endofunctor $g_{\Delta_M}$ of $F(M)$ canonically equivalent to the identity $\operatorname{id}_{F(M)}$.
(One should not expect the same brane structure on $\Delta_M$ to result in both functors $f_{\Delta_M}$ and~$g_{\Delta_M}$ being
the identity.)
\end{Remark}

\subsection{Closed cell correspondences}

Fix a~maximal critical point $p\in \mathfrak c$, so that the corresponding coisotropic cell is closed, and consider the
Hamiltonian reduction diagram
\begin{gather*}
\xymatrix{ M_p & \ar@{->>}[l]_-{q_p} C_p \ar@{^(->}[r]^-{i_p} & M,}
\end{gather*}
where $i_p$ is the inclusion of the coisotropic cell, and $q_p$ is the quotient by the integrable isotropic foliation determined
by $i_p$.

Def\/ine the {\em closed cell correspondence}
\begin{gather*}
\xymatrix{ \mathsf C_p \in F(M_p^{\operatorname{op}} \times M) }
\end{gather*}
to be the Lagrangian submanifold $C_p \subset M_p^{\operatorname{op}} \times M$ equipped with a~brane structure.
Since we are only sketching constructions, we will not specify the details of particular brane structures throughout what
follows.

We will sketch that the correspondence construction
\begin{gather*}
\xymatrix{ f_{\mathsf C_p}:\operatorname{Mod} F(M_p) \ar[r] & \operatorname{Mod} F(M) }
\end{gather*}
in fact should restrict to a~fully faithful representable functor
\begin{gather*}
\xymatrix{ \mathfrak i:F(M_p) \ar[r] & F(M).}
\end{gather*}

Observe that the geometric composition of a~Lagrangian brane $L\subset M_p$ with the correspondence $C_p$ is nothing more than
the pullback $i_p^{-1}L\subset M$.
The object $\mathfrak i (L)$ representing $f_{\mathsf C_p}(L)$ should be the pullback $i_p^{-1}L\subset M$ with an induced brane
structure.
This is a~geometric assertion whose proof should be similar to the statement that the diagonal brane gives the identity
correspondence.

To see why the functor should be fully faithful, consider the f\/iber $F_p = i_p^{-1}(p) \subset C_p$ and its cotangent bundle
$T^*F_p$.
We can f\/ind a~small open neighborhood ${\mathcal N}(C_p) \subset M$ of the coisotropic cell $C_p \subset M$ and a~symplectic
identif\/ication ${\mathcal N}(C_p) \simeq M_p \times U_p$, where the factor $U_p \subset T^*F_p$ is a~small open neighborhood of
the zero section.
Now observe that with this setup, the functor $\mathfrak i$ should be given by taking the product with the zero section
\begin{gather*}
\xymatrix{ F(M_p) \ar[r] & F(M_p \times U_p), & L \ar@{|->}[r] & L \times F_p.}
\end{gather*}
Following our perturbation framework, near inf\/inity we perturb such products by a~product of a~perturbation of $L$ and
a~perturbation of $F_p$.
For the latter, the perturbation will be to the graph of the dif\/ferential of a~function on $F_p$ which is linear at inf\/inity.
In particular, we can take the function on $F_p$ to have a~single critical point.

Now it should suf\/f\/ice to show that all disks involved in the calculation of morphisms between such products lie in ${\mathcal
N}(C_p)$.
Following for example~\cite{N, NZ}, this can be accomplished as follows.
First, by contracting the f\/irst factor with the Liouville vector f\/ield of~$M_p$, we can ensure arbitrarily small energy bounds
on the disks.
Then we can invoke~\cite{Sikorav} so that the energy bounds provide suf\/f\/icient diameter bounds.

\begin{Remark}
For a~f\/ixed marking $\Lambda_p\subset M_p$, consider the induced marking
\begin{gather*}
\Lambda_{p+} = K \cup q_p^{-1}(\Lambda_p) \subset M.
\end{gather*}

The functor $\mathfrak i$ should restrict to a~functor
\begin{gather*}
\xymatrix{ \mathfrak i: F_{\Lambda_p}(M_p) \ar[r] & F_{\Lambda_{p+}}(M).}
\end{gather*}
\end{Remark}

\begin{Remark}
One can alternatively apply the above discussion within the dual formalism of integral transforms suited to constructing left
adjoints.
As above, with the appropriate brane structure, the closed cell correspondence $\mathsf C_p \in F(M_p^{\operatorname{op}} \times M)$ should
give a~fully faithful functor~$g_{\mathsf C_p}$ that is also equivalent to~$\mathfrak i$.
(One should not expect the same brane structure on~$\mathsf C_p$ to lead to both functors $f_{\mathsf C_p}$ and $g_{\mathsf
C_p}$ being the functor~$\mathfrak i$.)
\end{Remark}

Next, there should exist a~right adjoint
\begin{gather*}
\xymatrix{ \mathfrak i^!:\operatorname{Mod} F(M) \ar[r] & \operatorname{Mod} F(M_p) }
\end{gather*}
given by the correspondence construction $f_{\mathsf C'_p}$ with an appropriate brane structure on $\mathsf C'_p$.
Once one knows that $\mathfrak i$ is given by the dual correspondence construction $g_{\mathsf C_p}$, it should be a~formal
consequence that its right adjoint is given by $f_{\mathsf C'_p}$.
Following patterns for constructible sheaves, we expect that $\mathfrak i^!$ should in fact be proper.

\begin{Remark}
One can alternatively apply the above discussion within the dual formalism of integral transforms suited to constructing left
adjoints.

First, one can regard the functor $\mathfrak i$ as a~proper continuous functor
\begin{gather*}
\xymatrix{ \operatorname{Mod} (F(M_p)^{\operatorname{op}}) \ar[r] & \operatorname{Mod} (F(M)^{\operatorname{op}}).}
\end{gather*}
Then there should exist an additional adjoint $\mathfrak i^*$ given by the correspondence construction $g_{\mathsf C'_p}$ with
an appropriate brane structure on $\mathsf C'_p$.
Once one knows that $\mathfrak i$ is given by the correspondence construction $f_{\mathsf C_p}$, it should be a~formal
consequence that the additional adjoint is given by $g_{\mathsf C'_p}$.

Moreover, following patterns for constructible sheaves, we expect that $\mathfrak i^*$ should in fact be proper.
Thus in turn it could be regarded as a~functor $\operatorname{Perf} F(M) \to \operatorname{Perf} F(M_p), $ which should provide
a~left adjoint to the original functor~$\mathfrak i$.
\end{Remark}

\subsection{Open complements}

We continue with a~maximal critical point $p\in \mathfrak c$, so that the corresponding coisotropic cell $C_p \subset M$ is
closed.
Now let us consider the open subset $\mathfrak s = \mathfrak c \setminus \{p\}$, and the corresponding open union of coisotropic
cells
\begin{gather*}
M_\mathfrak s = \coprod_{q \in \mathfrak s} C_q = M \setminus C_p.
\end{gather*}

Before proceeding further, it will be convenient to pin down the choice of a~def\/ining function $m_p:M \to [0, 1]$ for the
coisotropic cell $C_p \subset M$ more specif\/ically.
All of what follows can be achieved by applications of the Lagrangian and coisotropic neighborhood theorems~\cite{Gotay, Weina, Weinb} and the fact that $C_p\subset M$ is a~minimal unstable cell.

First, let us return to constructions seen in the preceding section.
Consider the f\/iber $F_p = i_p^{-1}(p) \subset C_p$ and its cotangent bundle $T^*F_p$.
We can f\/ind a~small open neighborhood ${\mathcal N}(C_p) \subset M$ and a~symplectic identif\/ication ${\mathcal N}(C_p) \simeq
M_p \times U_p$, where the factor $U_p \subset T^*F_p$ is a~small open neighborhood of the zero section.
Ref\/ining the cell $M_p$ of the f\/irst factor, we can choose a~Lagrangian $L_p \subset M_p$ and identif\/ications $L_p\simeq \mathbb
R^k$ and $T^*L_p \simeq M_p$, where $k = \dim M_p /2$.
Ref\/ining the normal second factor $U_p$, we can choose a~f\/iber $N_p \subset U_p$ of the projection $T^*F_p \to F_p$ and
identif\/ications $N_p \simeq \mathbb R^n$ and $U_p \simeq T^*N_p$, where $n = \dim M - \dim C_p$.
Thus altogether, we have identif\/ications ${\mathcal N}(C_p) \simeq T^*(L_p \times N_p) \simeq T^*(\mathbb R^{k} \times \mathbb
R^{n})$.
Finally, we can choose the def\/ining function $m_p:M \to [0, 1]$ so that outside of the open neighborhood ${\mathcal N}(C_p)
\subset M$, it is identically one, and near to $C_p \subset {\mathcal N}(C_p)$, it is simply the sum of the squares of the
coordinates of the normal factor $N_p \simeq \mathbb R^n$.

Consider the Hamiltonian function $\log m_{p}:M_{\mathfrak s} \to \mathbb R$, and the symplectomorphisms
\begin{gather*}
\xymatrix{ \Xi_{\mathfrak s!}:M_{\mathfrak s} \ar[r] & M_{\mathfrak s}, & \Xi_{\mathfrak s*}:M_{\mathfrak s} \ar[r] &
M_{\mathfrak s} }
\end{gather*}
resulting from the Hamiltonian f\/low of $\log m_{p}$ for negative unit time and unit time respectively.
To ensure good behavior near inf\/inity, it is technically useful to introduce a~suf\/f\/iciently small~\mbox{$\eta>0$} and the Hamiltonian
function $\log (m_{p} - \eta):M_{\mathfrak s, >\eta} \to \mathbb R$ on the domain $M_{\mathfrak s, >\eta} = \{x\in M_\mathfrak s
\, |\, m_p(x) >\eta\}$, and the symplectomorphisms
\begin{gather*}
\xymatrix{ \Xi_{\mathfrak s, \eta!}:M_{\mathfrak s, >\eta} \ar[r] & M_{\mathfrak s, >\eta}, & \Xi_{\mathfrak s,
\eta*}:M_{\mathfrak s, >\eta} \ar[r] & M_{\mathfrak s, >\eta} }
\end{gather*}
resulting from the Hamiltonian f\/low of $\log (m_{p} -\eta)$ for negative unit time and unit time respectively.
By subanalytic theory, any subanalytic subset $Y \subset \overline M$ must be transverse to the closed ball
$\overline{M_{\mathfrak s, \eta}} \subset \overline M$, given by the closure of the open ball $M_{\mathfrak s, \eta}= \{x\in
M_\mathfrak s \, |\, m_p(x) = \eta\}$, for suf\/f\/iciently small~$\eta>0$.

Now for f\/ixed $\dag = !$ or $*$, we expect that given a~brane $L\in F(M_\mathfrak s)$, for all suf\/f\/iciently small $\eta>0$, the
pushforwards
\begin{gather*}
  \Xi_{\mathfrak s, \eta \dag}(L \cap M_{\mathfrak s, \eta})
\end{gather*}
are well-def\/ined objects of $F(M)$ and all mutually dif\/feomorphic.

They also should all be dif\/feomorphic to the limit pushforward
\begin{gather*}
 \Xi_{\mathfrak s, 0^+ \dag}(L) = \lim\limits_{\eta\to 0^+} \Xi_{\mathfrak s, \eta \dag}(L \cap M_{\mathfrak s, \eta})
\end{gather*}
which also should be a~well-def\/ined object of $F(M)$.

We expect a~proof to proceed along the following lines.
By construction, $\Xi_{\mathfrak s, \eta \dag}(L \cap M_{\mathfrak s, \eta}) $ is dif\/feomorphic to $L$ and coincides with it on
a~common subset that is a~deformation retract of each.
Thus to see $\Xi_{\mathfrak s, \eta \dag}(L \cap M_{\mathfrak s, \eta}) $ is a~well-def\/ined brane, we should check that it is
closed in~$M$.
Suppose not: then there would be a~sequence of points along the boundary at inf\/inity of~$L$ at which $dm_p$ is a~negative scale
of the contact form $\alpha$ for the boundary of~$M_\mathfrak s$.
By the curve selection lemma, there should be a~curve $\gamma$ of such points, and hence $0> dm_p(\gamma') =c \alpha(\gamma')$
for some negative function $c$.
But this would contradict the fact that the boundary at inf\/inity of~$L$ must be isotropic.

Now we expect that for $\dag = !, *$, there is a~fully faithful functor
\begin{gather*}
\xymatrix{ F(M_\mathfrak s) \ar[r] & F(M) }
\end{gather*}
with continuous extension
\begin{gather*}
\xymatrix{ \mathfrak j_\dag:\operatorname{Mod} F(M_\mathfrak s) \ar[r] & \operatorname{Mod} F(M) }
\end{gather*}
given on Lagrangian branes $L\in F(M_\mathfrak s)$ by the limit
\begin{gather*}
\xymatrix{ \mathfrak j_ \dag L = \Xi_{\mathfrak s, 0^+ \dag}(L).
}
\end{gather*}

The existence of the functor $\mathfrak j_\dag$ and the fact that it is fully faithful should all follow from some simple
observations.
Consider a~f\/inite collection of branes in $M_\mathfrak s$ and f\/inite collection of $A_\infty$-compositions.
Observe that the construction $\Xi_{\mathfrak s, 0^+ \dag}$ is compatible with our perturbation framework: near inf\/inity,
perturbing forward in time is compatible with applying $\Xi_{\mathfrak s, 0^+ *}$, and backward in time with $\Xi_{\mathfrak s,
0^+ !}$.
Thus we can compare calculations before and after applying $\Xi_{\mathfrak s, 0^+ \dag}$ by collectively using the appropriate
direction of perturbation.
Next, observe that the relevant intersections and disks in all calculations are constrained to a~compact subset $K\subset
M_\mathfrak s$ thanks to diameter bounds depending only on a~smaller compact set.
Thus changes to the target and branes outside of $K\subset M_\mathfrak s$ are immaterial.
In particular, since the construction $\Xi_{\mathfrak s, 0^+ \dag}$ takes place outside of $M_{\mathfrak s, \eta}$, for
suf\/f\/iciently small $\eta$, it has no ef\/fect on any calculations.
Thus the extension of branes given by $\Xi_{\mathfrak s, 0^+ \dag}$ should be a~fully faithful embedding.

\begin{Remark}
Recall the induced markings
\begin{gather*}
\Lambda_{\mathfrak s} = \Lambda \cap M_\mathfrak s \subset M_\mathfrak s,
\qquad
\Lambda_p = q_p(\Lambda \cap C_p) \subset M_p,
\qquad
\Lambda_{p+} = \Lambda\cup q_p^{-1}(\Lambda_p) \subset M.
\end{gather*}

The functor $\mathfrak j_\dag$ should restrict to a~functor
\begin{gather*}
\xymatrix{ \mathfrak j_\dag:F_{\Lambda_\mathfrak s}(M_\mathfrak s) \ar[r] & F_{\Lambda_{p+}}(M).}
\end{gather*}
\end{Remark}

In the framework of Lagrangian correspondences, the functor $\mathfrak j_!$ should be given by the correspondence construction
$g_{\Gamma_{\mathfrak s !}}$ where we write $\Gamma_{\mathfrak s !} \subset M_{\mathfrak s}^{\operatorname{op}} \times M$ for the graph of
$\Xi_{\mathfrak s !}$ equipped with an appropriate brane structure.
This is a~geometric assertion whose proof should be similar to the statement that the diagonal brane gives the identity
correspondence.

Next, $\mathfrak j_!$ should admit a~right adjoint
\begin{gather*}
\xymatrix{ \mathfrak j^!:\operatorname{Mod} F(M) \ar[r] & \operatorname{Mod} F(M_\mathfrak s) }
\end{gather*}
given by the correspondence constructions $f_{\Gamma_{\mathfrak s !}'}$ with an appropriate brane structure on
$\Gamma_{\mathfrak s !}'$.
Once one knows that $\mathfrak j_!$ is given by the dual correspondence construction $g_{\Gamma_{\mathfrak s !}}$, it should be
a~formal consequence that its right adjoint is given by $f_{\Gamma_{\mathfrak s !}'}$.
Following patterns for constructible sheaves, we expect that $\mathfrak j^!$ should in fact be proper.

\begin{Remark}
One can alternatively apply the above discussion within the dual formalism of integral transforms suited to constructing left
adjoints.

First, the functor $\mathfrak j_*$ should be given by the correspondence construction $f_{\Gamma_{\mathfrak s *}}$ where we
write $\Gamma_{\mathfrak s *} \subset M_{\mathfrak s}^{\operatorname{op}} \times M$ for the graph of $\Xi_{\mathfrak s *}$ equipped with an
appropriate brane structure.

Then there should exist an additional adjoint $\mathfrak j^*$ given by the correspondence construc\-tion~$g_{\Gamma_{\mathfrak
s*}'}$ with an appropriate brane structure on $\Gamma_{\mathfrak s*}'$.
Moreover, following patterns for constructible sheaves, we expect that $\mathfrak j^*$ should in fact be proper.
Thus in turn it could be regarded as a~functor $\operatorname{Perf} F(M_\mathfrak s) \to \operatorname{Perf} F(M) $ which
should provide a~left adjoint to the original functor~$\mathfrak j_*$.
\end{Remark}

\subsection{Simple identities}

We mention here some simple relations that should hold between the functors introduced in the previous sections.
There should be vanishing
\begin{gather*}
\xymatrix{ \mathfrak j^!\mathfrak i \simeq 0:F(M_p) \ar[r] & F(M_\mathfrak s), & \mathfrak i^! \mathfrak j_*\simeq
0:F(M_\mathfrak s) \ar[r] & F(M_p) }
\end{gather*}
since thanks to our perturbation framework, given test branes, each of the above compositions should result in a~conf\/iguration
with non-intersecting branes.

\begin{Remark}
Within the dual formalism of integral transforms suited to constructing left adjoints, one should similarly have identities
\begin{gather*}
\mathfrak j^*\mathfrak i \simeq 0,
\qquad
\mathfrak i^* \mathfrak j_!\simeq 0.
\end{gather*}
\end{Remark}

\subsection{Geometry of diagonal}

We continue with the previous setup: a~maximal critical point $p\in \mathfrak c$ corresponding to a~closed coisotropic cell $C_p
\subset M$, and open subset $\mathfrak s = \mathfrak c \setminus \{p\}$ corresponding to the open complement
\begin{gather*}
M_\mathfrak s = \coprod_{q \in \mathfrak s} C_q = M \setminus C_p.
\end{gather*}

Let us sketch the existence of an open neighborhood ${\mathcal N}(\Delta_M \cup (C_p\times_{M_p} C_p)) \subset M^{\operatorname{op}} \times M$,
an open neighborhood ${\mathcal N}(T^*_M M \cup T^*_{C_p} M) \subset T^*M$, and a~symplectomorphism
\begin{gather*}
\xymatrix{ \psi:{\mathcal N}(\Delta_M \cup (C_p\times_{M_p} C_p)) \ar[r]^-\sim & {\mathcal N}(T^*_M M \cup T^*_{C_p} M),}
\end{gather*}
such that we have identif\/ications
\begin{gather*}
\psi(\Delta_M) = T^*_M M,
\qquad
\psi(C_p \times_{M_p} C_p) = T^*_{C_p} M,
\\
\psi(\Gamma_{\mathfrak s!} ) = -\Gamma_{d\log m_p},
\qquad
\psi(\Gamma_{\mathfrak s*} ) = \Gamma_{d\log m_p}.
\end{gather*}

For concreteness, let us return to the identif\/ications introduced earier.
Namely, we have an open neighborhood ${\mathcal N}(C_p) \subset M$, and identif\/ications
\begin{gather*}
{\mathcal N}(C_p) \simeq M_p \times U_p \simeq T^*(L_p \times N_p) \simeq T^*(\mathbb R^{k} \times \mathbb R^{n}),
\end{gather*}
where $L_p \subset M_p$, $N_p \subset U_p$ are Lagrangian cells.
Under these identif\/ications, the coisotropic cell~$C_p$ goes over to the coisotropic cell $T^*\mathbb R^k \times T^*_{\{0\}}
\mathbb R^n \simeq T^*(\mathbb R^k \times \mathbb R^n)|_{\mathbb R^k \times\{0\}}$.
Furthermore, the support of $d\log m_p$ is contained within ${\mathcal N}(C_p)$.

Let us simplify the notation and write $V= \mathbb R^{k} \times \mathbb R^{n}$, and $W = \mathbb R^{k} \times\{0\}$.
Let us further introduce notation for elements $x, y\in V$, $\xi, \eta \in V^*$, and consider the symplectomorphism
\begin{gather*}
\begin{split}
&\xymatrix{ (V\times V^*)^{\operatorname{op}} \times (V\times V^*) \simeq (T^*V)^{\operatorname{op}}
\times T^*V\ar[r]^-\sim & T^*(T^*V) \simeq (V \times V^*)
\times (V^* \times V),}
\\
& \xymatrix{ (x, \xi, y, \eta) \ar@{|->}[r] & \frac{1}{\sqrt{2}}(x+y, \xi+\eta, \eta - \xi, x-y ). }
\end{split}
\end{gather*}
It takes the coisotropic cell $(T^*V|_W) \times_{T^*W} (T^*V|_W)$ to the coisotropic cell $T^*_{T^*V|_W} (T^*V)$.

Now working f\/irst near $C_p \times_{M_p} C_p \subset M^{\operatorname{op}} \times M$, we can use the above constructions to def\/ine $\psi$
satisfying the requirements.
To extend it to a~neighborhood around the rest of $\Delta_M$, we can use the Lagrangian neighborhood theorem~\cite{Weina,
Weinb}.

Next def\/ine the {\em closed cell projector}
\begin{gather*}
  \Pi_p \in F(M^{\operatorname{op}} \times M)
\end{gather*}
to be the Lagrangian submanifold $C_p\times_{M_p} C_p \subset M^{\operatorname{op}} \times M$ equipped with an appropriate brane structure.

For $\dag = !, *$, def\/ine the {\em open complement projector}
\begin{gather*}
 \Pi_{\mathfrak s\dag} \in F(M^{\operatorname{op}} \times M)
\end{gather*}
to be the Lagrangian graph $\Gamma_{\Xi\dag} \subset M^{\operatorname{op}} \times M$ equipped with an appropriate brane structure.

Let us next sketch how inside of $F(M^{\operatorname{op}} \times M)$, we should have exact triangles of branes
\begin{gather*}
\xymatrix{ {\Pi_{p}} \ar[r] & \Delta_M \ar[r] & {\Pi_{\mathfrak s*}} \ar[r]^-{[1]} &,& {\Pi_{\mathfrak s!}} \ar[r] & \Delta_M
\ar[r] & {\Pi_p} \ar[r]^-{[1]} & }
\end{gather*}
with the appropriate brane structures.

On the one hand, given any open neighborhood ${\mathcal N}(\Delta_M) \subset M^{\operatorname{op}} \times M$, and any f\/inite number of
$A_\infty$-compositions among the branes in the assertion of the theorem, we can arrange so that all relevant geometry lies in
${\mathcal N}(\Delta_M)$.

On the other hand, using the microlocalization functor
\begin{gather*}
\xymatrix{ \mu_M:\operatorname{Sh}(M) \ar[r] & F(T^*M) }
\end{gather*}
we can transport the standard exact triangles
\begin{gather*}
\xymatrix{ i_!i^! {\mathcal D}_{M} \ar[r] & {\mathcal D}_M \ar[r] & j_* j^*{\mathcal D}_{M} \ar[r]^-{[1]} &, & j_!j^!k_{M} \ar[r]
& k_M \ar[r] & i_*i^*k_{M} \ar[r]^-{[1]} & },
\end{gather*}
where $k_M$ denotes the constant sheaf and ${\mathcal D}_M$ the Verdier dualizing sheaf, to exact triangles of branes
\begin{gather*}
\xymatrix{ {\Pi_{p}} \ar[r] & \Delta_M \ar[r] & {\Pi_{\mathfrak s*}} \ar[r]^-{[1]} &,& {\Pi_{\mathfrak s!}} \ar[r] & \Delta_M
\ar[r] & {\Pi_p} \ar[r]^-{[1]} & }
\end{gather*}
with the appropriate brane structures.

Furthermore, we can arrange so that all relevant geometry lies in any neighborhood ${\mathcal N}(\Delta_M \cup (C_p\times_{M_p}
C_p)) \subset M^{\operatorname{op}} \times M$.
Thus choosing a~neighborhood and symplectomorphism as sketched above, we obtain a~matching of all calculations.

It follows that any object $L \in \operatorname{Mod} F(M)$ should f\/it into an exact triangle
\begin{gather*}
\xymatrix{ \mathfrak i \mathfrak i^! L \ar[r] & L \ar[r] & \mathfrak j_* \mathfrak j^* L \ar[r]^-{[1]} & }.
\end{gather*}

To see this, apply the integral transform formalism to the f\/irst exact triangle of branes above to obtain an exact triangle
\begin{gather*}
\xymatrix{ f_{\Pi_{p}}(L) \ar[r] & L \ar[r] & f_{\Pi_{\mathfrak s*}}(L) \ar[r]^-{[1]} & }.
\end{gather*}
Then it remains to establish equivalences
\begin{gather*}
f_{\Pi_{p} } \simeq \mathfrak i \circ \mathfrak i^! \in \operatorname{Fun}(\operatorname{Mod} F(M), \operatorname{Mod} F(M)),
\\
f_{\Pi_{\mathfrak s*}} \simeq \mathfrak j_* \circ \mathfrak j^! \in \operatorname{Fun}(\operatorname{Mod} F(M),
\operatorname{Mod} F(M)).
\end{gather*}
These should admit direct verif\/ication or alternatively as a~non-compact variation on the theory of Lagrangian correspondences
of Wehrheim and Woodward~\cite{ww6, ww3,ww1, ww4, ww2, ww5}.

\begin{Remark}
The other expected triangle should result from taking the second triangle of branes above
\begin{gather*}
\xymatrix{ {\Pi_{\mathfrak s!}} \ar[r] & \Delta_M \ar[r] & {\Pi_p} \ar[r]^-{[1]} & & }
\end{gather*}
and applying the dual formalism of integral transforms suited to constructing left adjoints.
We expect all functors should be proper, and so for any object $L \in \operatorname{Mod} (F(M)^{\operatorname{op}})$, what should result is an
exact triangle
\begin{gather*}
\xymatrix{ \mathfrak j_! \mathfrak j^! L \ar[r] & L \ar[r] & \mathfrak i \mathfrak i^* L \ar[r]^-{[1]} & }
\end{gather*}
of objects of $F(M)$.
\end{Remark}

\subsection{A further identity}

We also expect to have an equivalence $\mathfrak j^!\simeq \mathfrak j^*$.
To see this, we seek a~natural equivalence
\begin{gather*}
\hom_{\operatorname{Perf} F(M_\mathfrak s)}(\mathfrak j^!L, P) \simeq \hom_{F(M)}(L, \mathfrak j_*P).
\end{gather*}
Since both sides vanish for objects of the form $L = \mathfrak i L'$, it suf\/f\/ices to assume that $L= \mathfrak j_!P'$, and to
establish a~natural equivalence
\begin{gather*}
\hom_{F(M_\mathfrak s)}(P', P) \simeq \hom_{F(M)}(\mathfrak j_! P', \mathfrak j_*P).
\end{gather*}
This should follow from the compatibility of the perturbation framework with the construction of the functors.

\section{Localization}

We return now to a~more traditional presentation.
The constructions of this section will not depend on the recollement pattern of diagrams~\eqref{eq functors} and~\eqref{eq
triangles} of the introduction as sketched in the preceding section.
But verifying their good properties in Section~\ref{sect compatibility} will.

We continue to suppose throughout this section that $(M, \theta, \eta, \sigma)$ is a~f\/ixed Weinstein target, so a~Weinstein
manifold $(M, \theta)$ with a~bicanonical trivialization $\eta$ and spin structure $\sigma$.
We will also suppose as well that each of its Weinstein cells $(M_p, \theta_p)$ is equipped with a~bicanonical trivialization
$\eta_p$ and the (necessarily) trivial spin structure $\sigma_p$.

\subsection{Preliminaries}

\begin{Definition}
The {\em conic topology} of $M$ is the category $M_{\rm con}$ with objects conic open subanalytic subsets of $M$ and morphisms
inclusions.
\end{Definition}

Recall that $\operatorname{st}_k$ denotes the $\infty$-category of small stable idempotent-complete $k$-linear $\infty$-categories.

\begin{Definition}
(1) A {\em $\operatorname{st}_k$-valued presheaf} on the conic topology of $M$ is an $\infty$-functor
\begin{gather*}
\xymatrix{ {\mathcal F}^{\rm pre}:M_{\rm con}^{\operatorname{op}} \ar[r] & \operatorname{st}_k.}
\end{gather*}

(2) A {\em $\operatorname{st}_k$-valued sheaf } on the conic topology of $M$ is a~continuous $\operatorname{st}_k$-valued presheaf.
\end{Definition}

\begin{Definition}
The {sheaf\/if\/ication} ${\mathcal F} = ({\mathcal F}^{\rm pre})^+$ of a~presheaf ${\mathcal F}^{\rm pre}$ is a~sheaf equipped with
a~universal presheaf morphism
\begin{gather*}
\xymatrix{ {\mathcal F}^{\rm pre} \ar[r] & {\mathcal F}.}
\end{gather*}
\end{Definition}

In our arguments, we will only need the following elementary functoriality.
Throughout what follows, ${\mathcal F}$ always denotes a~$\operatorname{st}_k$-valued presheaf or sheaf on the conic topology of $M$.

\begin{Definition}
Let $j:U\to M$ be the inclusion of a~conic open set.

The {\em restriction} $j^*{\mathcal F} = {\mathcal F}|_U$ is obtained by pullback along the induced functor $j:U_{\rm con} \to
M_{\rm con}$.
In other words, it assigns $ j^*{\mathcal F}(V) = {\mathcal F}(V) $ to any conic open set $V\subset U$.
\end{Definition}

\begin{Definition}
Let $\pi:M\to N$ be a~conic map.

The {\em pushforward} $\pi_*{\mathcal F}$ is obtained by pullback along the induced functor $\pi^{-1}:N_{\rm con} \to M_{\rm con}$.
In other words, it assigns $\pi_*{\mathcal F}(U) = {\mathcal F}(\pi^{-1}(U)) $ to any conic open set $U\subset N$.
\end{Definition}

\begin{Remark}
Restriction evidently commutes with sheaf\/if\/ication, while pushforward does not in general.
\end{Remark}

\begin{Remark}
Note that the global sections of a~presheaf or sheaf are simply its pushforward to a~point.
\end{Remark}

\begin{Definition}
The {\em support} of ${\mathcal F}$ is the smallest conic closed set $S\subset M$ such that ${\mathcal F}|_{M\setminus S}\simeq
0$.
\end{Definition}

\begin{Remark}
To see that the support is well-def\/ined, note that if $S_1, S_2\subset M$ are conic closed sets such that ${\mathcal
F}|_{M\setminus S_1} \simeq {\mathcal F}|_{M\setminus S_2}\simeq 0$, then the sheaf property implies ${\mathcal F}|_{M\setminus
(S_1 \cap S_2)} \simeq 0$.
\end{Remark}

\subsection{Construction of sheaf}

\begin{Definition}
Let $U\subset M$ be a~conic open subset.

We def\/ine the full subcategory of {\em $U$-null branes} $\operatorname{Null}(M, U) \subset \operatorname{Perf} F(M)$ to comprise
objects $L\in \operatorname{Perf} F(M)$ with singular support satisfying
\begin{gather*}
\xymatrix{ \operatorname{ss}(L) \cap U = \varnothing, }
\end{gather*}
or equivalently, null locus satisfying
\begin{gather*}
\xymatrix{ U \subset n(L).
}
\end{gather*}
\end{Definition}

\begin{Remark}
As a~special case of the above def\/inition, if we start with a~characteristic cone $\Lambda\subset M$, and then set $U=
M\setminus \Lambda$, we have by def\/inition $N(M, U) = \operatorname{Perf}_\Lambda F(M)$.
In particular, $\operatorname{Null}(M, \varnothing) = \operatorname{Perf} F(M)$.

Note as well in general $U_1\subset U_2$ implies $\operatorname{Null}(M, U_2) \subset \operatorname{Null}(M, U_1).
$\end{Remark}

\begin{Remark}
As a~consequence of our main results, we will deduce the nontrivial identity
\begin{gather*}
\xymatrix{ \operatorname{Null}(M, M) \simeq 0.
}
\end{gather*}
In other words, if a~brane $L\in \operatorname{Perf} F(M)$ has empty singular support $\operatorname{ss}(L) = \varnothing$, then it itself is
trivial $L\simeq 0$.
\end{Remark}

\begin{Definition}
\label{def sheaf}
(1) We def\/ine the $\operatorname{st}_k$-valued presheaf ${\mathcal F}_M^{\rm pre}$ on the conic topology of $M$ by the assignment
\begin{gather*}
\xymatrix{ {\mathcal F}_M^{\rm pre}(U) = \operatorname{Perf} F(M)/ \operatorname{Null}(M, U) }
\end{gather*}
for conic open subsets $U\subset M$.

(2) We def\/ine the $\operatorname{st}_k$-valued sheaf ${\mathcal F}_M$ of {\em localized branes} to be the sheaf\/if\/ication of ${\mathcal
F}_M^{\rm pre}$.
\end{Definition}

\begin{Remark}
In order to obtain the essential gluing property of a~sheaf, it is important that we follow Step (2) of the above def\/inition and
sheaf\/ify the naive quotient ${\mathcal F}_M^{\rm pre}$.
To see this explicitly, we recommend the reader continue on and consult Example~\ref{exr2} below where we illustrate this in the
easiest possible situation.
\end{Remark}

For conic open subsets $U\subset M$, we have the canonical localization morphism
\begin{gather*}
\xymatrix{ \operatorname{Loc}_U: \operatorname{Perf} F(M) \ar[r] & {\mathcal F}^{\rm pre}_M(U) \ar[r] & {\mathcal F}_M(U).}
\end{gather*}

\begin{Lemma}
\label{sheafss}
To each conic open subset $U\subset M$, and localized brane ${\mathcal L}\in {\mathcal F}_M(U)$, there is a~unique conic closed
subvariety $\operatorname{ss}_U({\mathcal L}) \subset U$ called the localized singular support characterized by the properties:
\begin{enumerate}\itemsep=0pt
\item[$(1)$] For conic open subsets $V\subset U\subset M$, we have compatibility with restriction
\begin{gather*}
\operatorname{ss}_V({\mathcal L}|_V) = \operatorname{ss}_U({\mathcal L}) \cap V.
\end{gather*}

\item[$(2)$] For a~brane $L \in \operatorname{Perf} F(M)$, we have compatibility with global singular support
\begin{gather*}
\operatorname{ss}_U(\operatorname{Loc}_U(L)) = \operatorname{ss}(L) \cap U.
\end{gather*}
\end{enumerate}
\end{Lemma}

\begin{proof}
The assertion is evident for sections of the presheaf ${\mathcal F}_M^{\rm pre}$.
Since conic closed subvarieties form a~sheaf, the assertion follows for sections of the sheaf ${\mathcal F}_M$.
\end{proof}

\subsection{Case of Weinstein cells}

Let $(N, \theta)$ be a~Weinstein cell.

\begin{Lemma}
If $L\in \operatorname{Perf} F(N)$ has empty singular support $\operatorname{ss}(L) = \varnothing$, then it itself is trivial $L\simeq 0$.
In other words, we have $\operatorname{Null}(N, N) \simeq 0$.
\end{Lemma}

\begin{proof}
By assumption, the unique zero $p \in N$ of the Liouville form lies in the null locus $n(L) \subset N$.
Hence by def\/inition, there exists a~conic open set $U\subset N$ containing $p$ such that for any test brane $P \in F(N)$ with
$P\subset U$, we have $\hom_{\operatorname{Perf} F(M)}(L, P) \simeq 0$.
But $N$ itself is the unique conic open set $U\subset N$ containing $p$.
Thus $\hom_{\operatorname{Perf} F(M)}(L, P) \simeq 0$ for any $P\in F(M)$, and hence $L \simeq 0$.
\end{proof}

\begin{Proposition}
Global localization is an equivalence
\begin{gather*}
\xymatrix{ \operatorname{Loc}_N:\operatorname{Perf} F(N) \ar[r]^-\sim & {\mathcal F}_N(N).}
\end{gather*}
\end{Proposition}

\begin{proof}
Note that $N$ itself is the unique conic open set containing the unique zero $p\in N$ of the Liouville form.
Hence any cover of $N$ by conic open sets must contain $N$ itself as a~constituent.
Thus the canonical map is an equivalence
\begin{gather*}
\xymatrix{ {\mathcal F}^{\rm pre}_N(N) \ar[r]^-\sim & {\mathcal F}_N(N).
}
\end{gather*}
Finally, by the previous lemma $\operatorname{Null}(N, N) \simeq 0$, hence the canonical map is an equivalence
\begin{gather*}
\xymatrix{ \operatorname{Perf} F(N) \ar[r]^-\sim &{\mathcal F}^{\rm pre}_N(N).} \tag*{\qed}
\end{gather*}
\renewcommand{\qed}{}
\end{proof}

\begin{Example}
\label{exr2}
Returning to Def\/inition~\ref{def sheaf}, let us see the impact of sheaf\/ifying via the easiest possible example.
Consider the two-dimensional Weinstein cell $M = \mathbb C$ with standard Liouville form $\theta$ and projectivization $M^\infty
\simeq S^1$.
Its core is the single point $K = \{0\} \subset \mathbb C$, and its ether is the complement $E = \mathbb C^* \subset \mathbb C$.

On the one hand, one can check that $\operatorname{Null}(\mathbb C, \mathbb C^*) \simeq 0$, and hence ${\mathcal F}_{M}^{\rm pre}(
\mathbb C^*) \simeq \operatorname{Perf} F(\mathbb C)$.
On the other hand, one can check that ${\mathcal F}_M(\mathbb C^*)\simeq \prod\limits_{x\in S^1} \operatorname{Perf} k$, where
the product is taken in $\operatorname{st}_k$.
The image of the canonical morphism
\begin{gather*}
\xymatrix{ \operatorname{Perf} F(\mathbb C) \simeq {\mathcal F}_{M}^{\rm pre}( \mathbb C^*)\ar[r] & {\mathcal F}_{M}( \mathbb C^*)
\simeq \prod\limits_{x\in S^1} \operatorname{Perf} k }
\end{gather*}
consists of sequences of objects whose underlying object (when we forget the labeling by points of $S^1$) is of the form
$V\oplus V[1] \in \operatorname{Perf} k$.
Informally speaking, the canonical morphism remembers the ``ends'' of branes near inf\/inity, and any exact Lagrangian curve in
$\mathbb C$ will have an even number of ends at inf\/inity.
More concretely, a~brane $L\in {\mathcal F}_{M}^{\rm pre}( \mathbb C^*) \simeq \operatorname{Perf} F(\mathbb C)$ supported along the
real line ${\mathbb R} \subset \mathbb C$, equipped with a~rank one local system, will have endomorphisms $k$, but its image in
$ {\mathcal F}_{M}( \mathbb C^*) \simeq \prod\limits_{x\in S^1} \operatorname{Perf} k$ will have endomorphisms $k\oplus k$.

Since ${\mathcal F}_{M}$ is a~sheaf, we see that ${\mathcal F}_{M}^{\rm pre}$ is not a~sheaf, and hence it was important that we
sheaf\/if\/ied it.
\end{Example}

\subsection{Compatibility with recollement}
\label{sect compatibility}

Now assume the recollement pattern of diagrams~\eqref{eq functors} and~\eqref{eq triangles} of the introduction as sketched in
the preceding section.

Let $(M, \theta)$ be a~Weinstein manifold.

Let $C\subset M$ be a~closed coisotropic cell and consider the Hamiltonian reduction diagram
\begin{gather*}
\xymatrix{ N & \ar@{->>}[l]_-{q} C \ar@{^(->}[r]^-{i} & M,}
\end{gather*}
where $i$ is the inclusion of the coisotropic cell, and $q$ is the quotient by the integrable isotropic foliation determined by
$i$.

Consider as well the complementary open
\begin{gather*}
\xymatrix{ j:M^\circ = M \setminus C \ar@{^(->}[r] & M.}
\end{gather*}

\begin{Proposition}
\label{nulldev}
Let $U^\circ\subset M^\circ \subset M$ be a~conic open subset.

Then the recollement functors restrict to a~diagram of adjunctions
\begin{gather*}
\xymatrix{ \ar[rr]^-{\mathfrak i} \operatorname{Perf} F(N) &&\operatorname{Null}(M, U^\circ)
\ar@/^2pc/[ll]_-{i^!}\ar@/_2pc/[ll]_-{\mathfrak i^*} \ar[rr]^-{\mathfrak j^! \simeq \mathfrak j^*} &&
\operatorname{Null}(M^\circ, U^\circ) \ar@/^2pc/[ll]_-{\mathfrak j_*} \ar@/_2pc/[ll]_-{\mathfrak j_!} }
\end{gather*}
\end{Proposition}

\begin{proof}
To avoid possible confusion, given a~conic subset $A\subset M^\circ \subset M$, we will write $A^\infty_M$ and~$A^\infty_{M^\circ} $ for its projectivizations as a~subset of $M$ and $M^\circ$ respectively.
Sinilarly, given a~closed subset $A\subset M^\circ \subset M$, we will write~$\partial_M^\infty A$ and
$\partial^\infty_{M^\circ} A$ for its boundaries at inf\/inity as a~subset of~$M$ and~$M^\circ$ respectively.

To see that the functor $\mathfrak i$ lands in $\operatorname{Null} (M, U^\circ)$, note f\/irst that $U^\circ \cap C = \varnothing$
and $(U^\circ)_M^\infty \cap C^\infty = \varnothing$.
Then any test brane $P\in F(M)$ with $P \subset U^\circ$ and $\partial_M^\infty P \subset (U^\circ)^\infty_M$ will not intersect
a~brane of the form $\mathfrak i( L) \in F(M)$, for any $L\in F(N)$.

From here, it suf\/f\/ices to see that the functor $\mathfrak j^*$ takes $\operatorname{Null} (M, U^\circ)$ to $\operatorname{Null}
(M^\circ, U^\circ)$.
Fix a~brane $L\in \operatorname{Null}(M, U^\circ)$, and any point $x\in U^\circ$.
We must conf\/irm that $x \in n(\mathfrak j^*(L))$.

Fix a~conic open set $V \subset U^\circ$ containing $x$ that exhibits $x\in n(L)$.
We will show that $V$ regarded as a~subset of $M^\circ$ also exhibits $x \in n(\mathfrak j^*(L))$.

Fix a~test brane $P\in F(M^\circ)$ with $P\subset V$ and $\partial_{M^\circ}^\infty V\subset V^\infty_{M^\circ}$.
Then we seek to show that
\begin{gather*}
\xymatrix{ \hom_{\operatorname{Perf} F(M^\circ)}(\mathfrak j^*(L), P)\simeq 0.
}
\end{gather*}
By adjunction, this is the same as to show that
\begin{gather*}
\xymatrix{ \hom_{\operatorname{Perf} F(M)}(L, \mathfrak j_*(P)) \simeq 0.
}
\end{gather*}
By construction, we have $\mathfrak j_*(P) \subset V$.
After the small perturbation required to compute the above morphisms, we also have $\partial^\infty_M(\mathfrak j_*(P)) \subset
V^\infty_M$.
Thus since $V$ exhibits $x\in n(L)$, the above morphisms vanish, and hence $V$ also exhibits $x\in n(\mathfrak j^*(L))$.
\end{proof}

\begin{Corollary}
\label{openequiv}

Restriction induces a~canonical equivalence
\begin{gather*}
\xymatrix{ \mathfrak j^! \simeq\mathfrak j^*: {\mathcal F}_{M}|_{M^\circ}\ar[r]^-\sim &{\mathcal F}_{M^\circ}.
}
\end{gather*}
\end{Corollary}

\begin{proof}
It suf\/f\/ices to show the analogous statement for presheaves
\begin{gather*}
\xymatrix{ {\mathcal F}^{\rm pre}_{M}|_{M^\circ}\ar[r]^-\sim &{\mathcal F}^{\rm pre}_{M^\circ}.
}
\end{gather*}
In other words, for conic open subsets $U^\circ\subset M^\circ$, it suf\/f\/ices to show that the restriction descends to compatible
equivalences
\begin{gather*}
\xymatrix{ \operatorname{Perf} F(M)/ \operatorname{Null}(M, U^\circ) \ar[r]^-\sim & \operatorname{Perf} F(M^\circ)/
\operatorname{Null}(M^\circ, U^\circ).}
\end{gather*}
This follows immediately from the recollement compatibility of Proposition~\ref{nulldev}.
\end{proof}

\begin{Proposition}
\label{nulldev2}
Let $U\subset M$ be a~conic open subset containing $C$, and let $U^\circ\subset M^\circ$ denote the conic open subset $U^\circ =
U \cap M^\circ$.

Then the recollement functors restrict to a~diagram of equivalences
\begin{gather*}
\xymatrix{ \operatorname{Null}(M, U) \ar[rr]_-\sim^-{\mathfrak j^! \simeq \mathfrak j^*} && \operatorname{Null}(M^\circ,
U^\circ) \ar@/^2pc/[ll]_-{\mathfrak j_*}^-\sim \ar@/_2pc/[ll]_-{\mathfrak j_!}^-\sim }
\end{gather*}
\end{Proposition}

\begin{proof}
It is evident that $\mathfrak j_!, \mathfrak j_*$ take $\operatorname{Null} (M^\circ, U^\circ)$ to $\operatorname{Null}(M, U)$,
and also that $\operatorname{Null} (M, U) \cap \mathfrak i(\operatorname{Perf} F(M)) = 0$.
\end{proof}

\begin{Corollary}
\label{sheafadj}
Let $U\subset M$ be a~conic open subset containing $C$, and let $U^\circ\subset M^\circ$ denote the conic open subset $U^\circ =
U \cap M^\circ$.

Then the recollement functors induce a~diagram of adjunctions
\begin{gather*}
\xymatrix{ \ar[rr]^-{\mathfrak i_! \simeq \mathfrak i_*} \operatorname{Perf} F(N) &&{\mathcal F}^{\rm pre}_{M}(U)
\ar@/^2pc/[ll]_-{i^!}\ar@/_2pc/[ll]_-{\mathfrak i^*} \ar[rr]^-{\mathfrak j^! \simeq \mathfrak j^*} &&{\mathcal
F}^{\rm pre}_{M^\circ}(U^\circ) \ar@/^2pc/[ll]_-{\mathfrak j_*} \ar@/_2pc/[ll]_-{\mathfrak j_!} }
\end{gather*}

Consequently, we have an equivalence
\begin{gather*}
{\mathcal F}^{\rm pre}_M(U) \simeq \operatorname{Mod} _T\big({\mathcal F}^{\rm pre}_{M^\circ}(U^\circ) \oplus \operatorname{Perf}F(N)\big),
\end{gather*}
where $T = RL\in \operatorname{End}({\mathcal F}_{M^\circ}^{\rm pre}(U) \oplus \operatorname{Perf} F(N))$ is the monad of the
adjunction
\begin{gather*}
\xymatrix{ L =\mathfrak j_! \oplus\mathfrak i_!: {\mathcal F}_{M^\circ}^{\rm pre}(U^\circ) \oplus \operatorname{Perf}
F(N)\ar@<0.5ex>[r] &\ar@<0.5ex>[l] {\mathcal F}_M^{\rm pre}(U):R =\mathfrak j^! \oplus \mathfrak i^!.}
\end{gather*}
\end{Corollary}

\begin{proof}\sloppy 
The f\/irst part follows immediately from the recollement compatibility of Proposi\-tion~\ref{nulldev2}.
The second part is an immediate application of the Barr--Beck theorem.
\end{proof}

\begin{Remark}
While the abstract monadic language is convenient, little of the sophisticated theory it represents is needed.
More simply, we can say that ${\mathcal F}^{\rm pre}_M(U)$ is equivalent to the $\infty$-category of triples $L^\circ \in {\mathcal
F}^{\rm pre}_{M^\circ}(U^\circ)$, $L_N \in \operatorname{Perf} F(N)$ and a~morphism
\begin{gather*}
r \in \hom_{\operatorname{Perf} F(N)}(\mathfrak i^! \mathfrak j_!L^\circ, L_N).
\end{gather*}
\end{Remark}

\subsection{Global sections}

Continue to assume the recollement pattern of diagrams~\eqref{eq functors} and~\eqref{eq triangles} of the introduction.

Let $(M, \theta)$ be a~Weinstein manifold.

Let $C\subset M$ be a~closed coisotropic cell and consider the Hamiltonian reduction diagram
\begin{gather*}
\xymatrix{ N & \ar@{->>}[l]_-{q} C \ar@{^(->}[r]^-{i} & M,}
\end{gather*}
where $i$ is the inclusion of the coisotropic cell, and $q$ is the quotient by the integrable isotropic foliation determined by
$i$.

Consider as well the complementary open
\begin{gather*}
\xymatrix{ j:M^\circ = M \setminus C \ar@{^(->}[r] & M.}
\end{gather*}

To f\/ind a~natural context for Corollaries~\ref{openequiv} and~\ref{sheafadj}, let us consider the conic quotient
\begin{gather*}
\xymatrix{ \pi: M \ar[r] & M^\sim = M^\circ \cup *,}
\end{gather*}
where we collapse $C\subset M$ to a~point denoted by $*$.
Observe that the inverse-image under $\pi$ provides an equivalence from the category of conic open sets $U^\sim \subset M^\sim$
to the category of conic open sets $U\subset M$ such that $U \cap C$ is either all of $C$ or empty.

Let us introduce the pushforward presheaves $\pi_*{\mathcal F}^{\rm pre}_M$ and $\pi_*j_*{\mathcal F}_{M^\circ}^{\rm pre}$, and denote
by $\operatorname{Perf} F(N)_*$ the skyscraper sheaf with f\/iber $\operatorname{Perf} F(N)$ supported at $*$.
Then we can reformulate Corolla\-ries~\ref{openequiv} and~\ref{sheafadj} as a~diagram of adjunctions of presheaves
\begin{gather*}
\xymatrix{ \ar[rr]^-{\mathfrak i_! \simeq \mathfrak i_*} \operatorname{Perf} F(N)_* && \pi_*{\mathcal F}^{\rm pre}_{M}
\ar@/^2pc/[ll]_-{i^!}\ar@/_2pc/[ll]_-{\mathfrak i^*} \ar[rr]^-{\mathfrak j^! \simeq \mathfrak j^*} &&\pi_* j_*{\mathcal
F}^{\rm pre}_{M^\circ} \ar@/^2pc/[ll]_-{\mathfrak j_*} \ar@/_2pc/[ll]_-{\mathfrak j_!} }
\end{gather*}

Consequently, we have an equivalence
\begin{gather*}
 \pi_*{\mathcal F}^{\rm pre}_M \simeq \operatorname{Mod} _T\big(\pi_* j_* {\mathcal F}^{\rm pre}_{M^\circ} \oplus
\operatorname{Perf} F(N)_*\big),
\end{gather*}
where $T = RL\in \operatorname{End}(\pi_* j_* {\mathcal F}_{M^\circ}^{\rm pre} \oplus \operatorname{Perf} F(N)_*)$ is the monad of
the adjunction
\begin{gather*}
\xymatrix{ L =\mathfrak j_! \oplus\mathfrak i_!: \pi_* j_* {\mathcal F}_{M^\circ}^{\rm pre} \oplus \operatorname{Perf}
F(N)_*\ar@<0.5ex>[r] &\ar@<0.5ex>[l] \pi_*{\mathcal F}_M^{\rm pre}:R =\mathfrak j^! \oplus \mathfrak i^!.}
\end{gather*}

In concrete terms, to any conic open set $U^\sim \subset M^\sim$, we have that $\pi_*{\mathcal F}^{\rm pre}_M(U^\sim)$ is equivalent
to the $\infty$-category of triples $L^\circ \in {\mathcal F}^{\rm pre}_{M^\circ}(U^\circ)$ where $U^\circ = U^\sim \cap M^\circ$,
$L_N \in \operatorname{Perf} F(N)$ nonzero only if $*\in U^\sim$, and a~morphism
\begin{gather*}
r \in \hom_{\operatorname{Perf} F(N)}\big(\mathfrak i^! \mathfrak j_!L^\circ, L_N\big).
\end{gather*}

Now we will check that the above description is compatible with sheaf\/if\/ication.

\begin{Lemma}
The canonical morphism is an equivalence of sheaves
\begin{gather*}
\xymatrix{ (\pi_*{\mathcal F}_M^{\rm pre})^+ \ar[r]^-\sim & \pi_*{\mathcal F}_M.}
\end{gather*}
\end{Lemma}

\begin{proof}
This is evident over the open set $M^\circ \subset M^\sim$.
It suf\/f\/ices to check that the canonical morphism induces an equivalence on stalks at $*$.
This follows from the further observation that any conic open set $U\subset M$ containing the zero $p\in \pi^{-1}(*) = C$ in
fact contains all of $C$.
\end{proof}

\begin{Theorem}
\label{sheafmonad}
The pushforward $\pi_*{\mathcal F}_M$ admits the canonical description
\begin{gather*}
\pi_*{\mathcal F}_M \simeq \operatorname{Mod} _T(\pi_* j_*{\mathcal F}_{M^\circ} \oplus \operatorname{Perf} F(N)_*),
\end{gather*}
where $T = RL\in \operatorname{End}(\pi_* j_*{\mathcal F}_{M^\circ} \oplus \operatorname{Perf} F(N)_*)$ is the monad of the
adjunction
\begin{gather*}
\xymatrix{ L =\mathfrak j_! \oplus\mathfrak i_!:\pi_* j_* {\mathcal F}_{M^\circ}^{\rm pre} \oplus \operatorname{Perf}
F(N)_*\ar@<0.5ex>[r] &\ar@<0.5ex>[l] {\mathcal F}_M:R =\mathfrak j^! \oplus \mathfrak i^!.}
\end{gather*}
\end{Theorem}

\begin{proof}
By our previous results reformulated above, the sheaf\/if\/ication $(\pi_*{\mathcal F}_M^{\rm pre})^+$ clearly admits the asserted
description.
Thus by the previous lemma, the pushforward $\pi_*{\mathcal F}_M$ does as well.
\end{proof}

\begin{Corollary}
Global localization is an equivalence
\begin{gather*}
\xymatrix{ \operatorname{Loc}_M: \operatorname{Perf} F(M) \ar[r]^-\sim & {\mathcal F}_M(M).}
\end{gather*}
\end{Corollary}

\begin{proof}
Note that $\pi_*{\mathcal F}_M(M^\sim) \simeq {\mathcal F}_M(M)$ and $\pi_* j_*{\mathcal F}_{M^\circ}(M^\sim) \simeq {\mathcal
F}_{M^\circ}(M^\circ)$.
By induction, global localization is an equivalence on the open Weinstein manifold
\begin{gather*}
\xymatrix{ \operatorname{Loc}_{M^\circ}: \operatorname{Perf} F(M^\circ) \ar[r]^-\sim & {\mathcal F}_{M^\circ}(M^\circ).}
\end{gather*}
Hence by Theorem~\ref{sheafmonad}, we have an equivalence on global sections
\begin{gather*}
\xymatrix{ {\mathcal F}_M (M) \simeq \operatorname{Mod} _T(\operatorname{Perf} F(M^\circ) \oplus \operatorname{Perf} F(N)),}
\end{gather*}
where $T = RL\in \operatorname{End}(\operatorname{Perf} F(M^\circ) \oplus \operatorname{Perf} F(N))$ is the monad of the
adjunction
\begin{gather*}
\xymatrix{ L =\mathfrak j_! \oplus\mathfrak i_!:\operatorname{Perf} F(M^\circ) \oplus \operatorname{Perf} F(N)\ar@<0.5ex>[r]
&\ar@<0.5ex>[l] {\mathcal F}_M(M):R =\mathfrak j^! \oplus \mathfrak i^!.}
\end{gather*}
Comparison with the similar monadic description of $\operatorname{Perf} F(M)$ yields the theorem.
\end{proof}

{\samepage \begin{Corollary}
For $L\in \operatorname{Perf} F(M)$, if $\operatorname{ss}(L) = \varnothing$, then $L\simeq 0$.
\end{Corollary}

\begin{proof}
The localization of $L$ is a~null brane for any conic open set.
\end{proof}}

\begin{Remark}
We will not need the following discussion but include it to help further orient the interested reader.

One might ask whether a~recollement description similar to that of Theorem~\ref{sheafmonad} might exist for the sheaf ${\mathcal
F}_M$ itself rather than its pushforward $\pi_*{\mathcal F}_M$.
The immediate answer is negative since the key functors $\mathfrak j_!, \mathfrak j_*$ are not local.
But one need not pass all the way to the quotient $M\to M^\sim$ induced by the collapse $C\to *$.
Rather it is possible to pass to the intermediate quotient $M \to{}'M^\sim$ induced by the natural collapse $C\to N$.
More broadly, the sheaf which admits a~natural recollement pattern is the pushforward of ${\mathcal F}_M$ along the quotient of
$M$ where each coisotropic cell is collapsed to its corresponding Weinstein cell.

In another direction, one might also ask which aspects of the recollement pattern can be lifted to the sheaf ${\mathcal F}_M$
itself.
First, we can consider the full subsheaf $i_!q^*{\mathcal F}_{N} \subset {\mathcal F}_M$ generated by objects of the form
$\mathfrak i(L)\in \operatorname{Perf} F(M)$, for objects $L\in \operatorname{Perf} F(N)$.
There is a~canonical equivalence on global sections
\begin{gather*}
\xymatrix{ \operatorname{Perf} F(N) \simeq {\mathcal F}_{N}(N) \ar[r]^-\sim & i_!q^*{\mathcal F}_{N}(M).}
\end{gather*}

Second, we have seen that there is a~canonical morphism ${\mathcal F}_M \to j_* {\mathcal F}_{M^\circ}$ that induces an
equivalence
\begin{gather*}
\xymatrix{ {\mathcal F}_{M}|_{M^\circ}\ar[r]^-\sim &{\mathcal F}_{M^\circ}.
}
\end{gather*}
Unfortunately, as mentioned above, there are no evident adjoint maps of sheaves.
The following related construction provides a~partial solution.
For $\dag = !$ or $*$, we can consider the full subsheaf ${\mathcal F}_{M^\circ \dag} \subset {\mathcal F}_M$ generated by
objects of the form $\mathfrak j_\dag L\in \operatorname{Perf} F(M)$, for objects $L\in \operatorname{Perf} F(M^\circ)$.
(We caution the reader that ${\mathcal F}_{M^\circ \dag} $ is not the same as the pushforward $j_\dag {\mathcal F}_{M}$.) Then
restriction induces a~canonical equivalence
\begin{gather*}
\xymatrix{ {\mathcal F}_{M\dag }|_{M^\circ}\ar[r]^-\sim &{\mathcal F}_{M^\circ}.
}
\end{gather*}
Furthermore, it also induces a~canonical equivalence on global sections
\begin{gather*}
\xymatrix{ {\mathcal F}_{M\dag }(M) \ar[r]^-\sim & {\mathcal F}_{M^\circ}(M^\circ).
}
\end{gather*}
\end{Remark}

\subsection{Prescribed support} Continue to assume the recollement pattern of diagrams~\eqref{eq functors} and~\eqref{eq
triangles} of the introduction.

Now let us introduce a~characteristic cone $\Lambda \subset M$ so that we have a~marked Weinstein manifold $(M, \theta,
\Lambda)$.
Recall by Lemma~\ref{sheafss}, we have the notion of singular support for localized branes.

\begin{Definition}

We def\/ine the full subsheaf ${\mathcal F}_\Lambda \subset {\mathcal F}_M$ to consist of those localized branes ${\mathcal L}
\subset {\mathcal F}_M(U)$ such that
\begin{gather*}
 \operatorname{ss}_U({\mathcal L}) \subset \Lambda \cap U
\end{gather*}
for any conic open set $U\subset M$.
\end{Definition}

We now have the following assertion from the introduction.

\begin{Theorem}
Assume the recollement pattern of diagrams~\eqref{eq functors} and~\eqref{eq triangles} of the introduction.

Let $(M, \theta, \Lambda)$ be a~marked Weinstein manifold.

The ${st}_k$-valued sheaf ${\mathcal F}_\Lambda$ on the conic topology of $M$ has the following properties:
\begin{enumerate}\itemsep=0pt
\item[$(1)$] The support of ${\mathcal F}_\Lambda$ is the characteristic cone $\Lambda \subset M$.

\item[$(2)$] The global sections of ${\mathcal F}_\Lambda$ are canonically equivalent to $\operatorname{Perf}_\Lambda(M)$.

\item[$(3)$] The restriction of ${\mathcal F}_\Lambda$ to an open Weinstein submanifold $M^\circ \subset M$ is canonically equivalent to
the sheaf ${\mathcal F}_{\Lambda^\circ}$ constructed with respect to $\Lambda^\circ = \Lambda \cap M^\circ$.

\item[$(4)$] For each zero $p\in \mathfrak c$, the sections of ${\mathcal F}_\Lambda$ lying strictly above the unstable cell $C_p \subset
M$ are canonically equivalent to $\operatorname{Perf}_{\Lambda_p} (M_p)$.
\end{enumerate}
\end{Theorem}

\begin{Example}
Let us return to the setting of Example~\ref{exr2} and add a~characteristic cone to the mix.

Recall the two-dimensional Weinstein cell $M = \mathbb C$ with standard Liouville form $\theta$ and projectivization $M^\infty
\simeq S^1$.
Its core is the single point $K = \{0\} \subset \mathbb C$, and its ether is the complement $E = \mathbb C^* \subset \mathbb C$.
Any characteristic cone $\Lambda \subset \mathbb C$ will be the union of $K = \{0\}$ with f\/initely many rays.
Let us denote by $\Lambda_n \subset \mathbb C$ the characteristic cone with $n$ rays, for $n = 0, 1, 2, \ldots$.

Then as is well known, $\operatorname{Perf}_{\Lambda_n} F(\mathbb C)$ is equivalent to f\/inite-dimensional modules over the
$A_{n-1}$-quiver (in particular, for $n=0$ and $n=1$, it is the zero category).
This is also the stalk of the sheaf ${\mathcal F}_{\Lambda_n}$ at the point $0\in \mathbb C$.
Its stalk at other points $x\in \mathbb C$ is (not necessarily canonically) equivalent to $\operatorname{Perf} k$ when $x\in
\Lambda_n$, and the zero category otherwise.
\end{Example}

\subsection*{Acknowledgements} I am indebted to D.~Ben-Zvi, P.~Seidel and E.~Zaslow for the impact they have had on my thinking
about symplectic and homotopical geometry.
I am grateful to T.~Perutz and D.~Treumann for many stimulating discussions, both of a~technical and philosophical nature.
I am grateful to M.~Abouzaid and D.~Auroux for their patient explanations of foundational issues and related questions in mirror
symmetry.
I~am also grateful to the anonymous referees for their thoughtful reading and generous investment in improving the paper.
I~would like to thank A.~Preygel for sharing his perspective on ind-coherent sheaves.
I~am also pleased to acknowledge the moti\-vating inf\/luence of a~question asked by C.~Teleman at ESI in Vienna in January 2007.
Finally, I~am grateful to the participants of the June 2011 MIT RTG Geometry retreat for their inspiring interest in this topic.

This work was supported by NSF grant DMS-0600909.

\pdfbookmark[1]{References}{ref}
\LastPageEnding

\end{document}